\documentclass[12pt]{article}

\usepackage[english]{babel}
\usepackage{amsmath,amssymb,amsthm}
\usepackage[mathscr]{eucal}
\usepackage{bbm}

\newcommand{\C}{\mathbb{C}}
\newcommand{\D}{\mathbb{D}}
\newcommand{\Disk}{\mathbb{D}}
\newcommand{\E}{\mathbb{E}}
\renewcommand{\H}{\mathbb{H}}
\newcommand {\Half} {\mathbb{H}}
\newcommand{\N}{\mathbb{N}}
\renewcommand{\P}{\mathbb{P}}
\newcommand{\R}{\mathbb{R}}

\newcommand{\Z}{\mathbb{Z}}

\newcommand{\A}{\mathcal{A}}
\newcommand{\B}{\mathcal{B}}

\newcommand{\Dr}[1]{\mathcal{D}^{#1}}
\newcommand{\scd}{\mathcal{D}}

\newcommand{\F}{\mathcal{F}}
\newcommand{\K}{\mathcal{K}}
\newcommand{\curves}{\mathcal{K}}
\newcommand{\Exc}{\mathcal{K}}
\newcommand{\M}{\mathcal{M}}
\newcommand{\finmeas}{\mathcal{M}}
\newcommand{\probmeas}{\mathcal{PM}}

\newcommand{\Square}{\mathcal{S}}

\newcommand{\BX}{\mathcal{X}}

\newcommand{\w}{\omega}
\newcommand{\e}{\varepsilon}
\newcommand{\eps}{\varepsilon}
\renewcommand{\epsilon}{\varepsilon}
\renewcommand{\phi}{\varphi}

\newcommand{\X}{\Xi}
\newcommand{\dn}{\delta_n}
\newcommand{\dN}{\delta_N}

\newcommand{\G}{\Gamma}
\newcommand{\U}{\Upsilon}
\newcommand{\GD}{\Gamma_{\D}}
\newcommand{\UD}{\Upsilon_{\D}}
\newcommand{\GN}{\Gamma_N}
\newcommand{\UN}{\Upsilon_N}
\newcommand{\tGN}{\tilde{\Gamma}_N}
\newcommand{\tUN}{\tilde{\Upsilon}_N}
\newcommand{\trunc}[2]{\Theta_{#1}^{#2}\gamma}
\newcommand{\truncopr}[2]{\Theta_{#1}^{#2}}

\renewcommand{\d}{\mathrm{d}}
\newcommand{\n}{\mathbf{n}}
\newcommand{\RW}{\mathsf{rw}}
\newcommand{\lenD}[1]{\ell_{#1}}
\newcommand{\brscale}{\Psi}

\newcommand{\rad}{\operatorname{rad}}
\newcommand{\inrad}{\operatorname{inrad}}
\newcommand{\dist}{\operatorname{dist}}
\newcommand{\diam}{\operatorname{diam}}
\newcommand{\cl}{\operatorname{cl}}
\newcommand{\interior}{\operatorname{int}}
\newcommand{\abs}{\operatorname{abs}}
\newcommand{\osc}{\operatorname{osc}}
\newcommand{\sep}{\operatorname{sep}}
\newcommand{\spr}{\operatorname{spr}}

\newcommand{\metric}{\mathbbm{d}}
\newcommand{\metricnoparam}{d^{*}_{\K}}
\newcommand{\metricparam}{d_{\K}}
\newcommand{\prohorov}{\wp}
\newcommand{\param}{\overset{\mathsf{par}}\sim}
\newcommand{\cara}{\overset{\mathsf{cara}}\to}

\newcommand{\bd}{\partial}  
\newcommand{\p}{\partial}

\DeclareMathAlphabet{\mathpzc}{OT1}{pzc}{m}{it}

\renewcommand{\Im}{\mathpzc{Im}}

\newcommand{\defines}[1]{\textit{#1}\index{#1}}

\newtheorem{theorem}{Theorem}[section]
\newtheorem{proposition}[theorem]{Proposition}
\newtheorem{lemma}[theorem]{Lemma}
\newtheorem{corollary}[theorem]{Corollary}
\newtheorem{importantremark}[theorem]{Important Remark}

\theoremstyle{definition}
\newtheorem{definition}[theorem]{Definition}
\newtheorem{example}[theorem]{Example}

\title{On the scaling limit of simple random walk excursion measure in the plane}

\date{}

\author{
Michael J.~Kozdron\footnote{Research supported in part by the Natural Sciences and Engineering Research Council of Canada.}\\
{\small Department of Mathematics \& Statistics, College West 307.14}\\
{\small University of Regina, Regina, SK S4S 0A2 Canada}\\ 
{\small \texttt{kozdron@math.uregina.ca}}\\        
}

\begin{document}

\maketitle

\begin{abstract}
The Brownian excursion measure is a conformally invariant infinite measure on curves. It figured prominently in one of the first major applications of SLE, namely the explicit calculations of the planar Brownian intersection exponents from which the Hausdorff dimension of the frontier of the Brownian path could be computed (Lawler, Schramm, and Werner, 2001). In this paper we define the simple random walk excursion measure and show that for any bounded, simply connected Jordan domain $D$, the simple random walk excursion measure on $D$ converges in the scaling limit to the Brownian excursion measure on $D$.
\end{abstract}



\newpage

\section{Introduction}\label{Sect1}

A number of mathematically simplistic lattice models, including the self-avoiding random walk~\cite{flory}, have been introduced in an attempt to better understand critical phenomena in two-dimensional statistical physics. While these models have been studied for several decades, little progress had been made until recently. The introduction of the Schramm-Loewner evolution, a new family of conformally invariant distributions on random curves, has led to a plethora of exciting results about the scaling limits of these models at criticality. For example, the scaling limits of loop-erased random walk~\cite{LawSW9,Zhan}, uniform spanning trees~\cite{LawSW9}, and site percolation on the triangular lattice~\cite{camia,Smirnov} can now be described using SLE.
 
One of the first successes, however, of the SLE program was the determination of the intersection exponents for random walk and Brownian motion, and the establishment of Mandelbrot's conjecture that the Hausdorff dimension of the frontier of the planar Brownian path is $4/3$. (See~\cite{LSW4} for a survey of this work.) The Brownian excursion measure, a conformally invariant infinite measure on curves which had been introduced in previous work by Lawler and Werner~\cite{LawW1}, figured prominently in the explicit calculations of the intersection exponents.

The goal of this present paper is to construct a discrete object, the simple random walk excursion measure, which has the Brownian excursion measure as its scaling limit. Of course, the convergence of simple random walk on $\Z^2$ to Brownian motion in $\C$ has been known since Donsker's theorem of 1951.  However, what had not been established was a strong version of this result which holds for random walk and Brownian motion on any simply connected domain where the errors do not depend on the smoothness of the boundary. By proving in the present paper that for any bounded, simply connected Jordan domain, the scaling limit of discrete excursion measure is Brownian excursion measure, we establish such a result. 

\subsection{Main results}

We begin with a discussion of the main results, leaving some of the precise statements to later sections. Our concern will be exclusively two dimensional, so we will identify $\C \cong \R^2$ in the usual way, and write any of $w$, $x$, $y$, or $z$ for points in $\C$. A \defines{domain} $D \subset \C$ is an open and connected set; write $\D := \{z \in \C : |z|< 1\}$ for the open unit disk, and $\H := \{z \in \C : \Im(z)>0\}$ for the upper half plane. A standard complex Brownian motion will be denoted $B_t$, $t\ge 0$, and $S_n$, $n=0, 1, \ldots$, will denote two-dimensional simple random walk, both started at the origin unless otherwise noted. We will generally use $T$ for stopping times for Brownian motion and $\tau$ for stopping times for random walk, and write $\E^x$ and $\P^x$ for expectations and probabilities, respectively, assuming $B_0=x$ or $S_0=x$. 

A subset $A \subset \Z^2$ is said to be \defines{simply connected} if both $A$ and $\Z^2 \setminus A$ are non-empty and connected. Write the \defines{(outer) boundary} of $A$ as $\bd A := \{ z \in \Z^2 \setminus A : \dist(z, A)=1 \}$. An \defines{excursion in $A$} is  a path $\w := [\w_0, \w_1, \ldots, \w_k]$ with $|\w_j - \w_{j-1}| = 1$ for all $j$; $\w_0$, $\w_k \in \bd A$; and $w_1, \ldots, \w_{k-1} \in A$. It is implicit that $2 \le k < \infty$; the \defines{length} of $\w$ is $|\w| := k$. We can view excursions of length $k$ as curves $\w:[0, k] \to \C$ by linear interpolation. Write $\Exc_A$ for the set of excursions in $A$, and define the \defines{simple random walk excursion measure} as the measure on $\Exc_A$ which assigns measure $4^{-k}$ to each length $k$ excursion in $A$. That is, the excursion measure of $\w = [\w_0, \w_1, \ldots, \w_k]$ is the probability that the first $k$ steps of a simple random walk starting at $\w_0$ are the same as $\w$. Let $D \subset \C$ be a bounded simply connected domain containing the origin, and for each $N < \infty$, let $D_N$ denote the connected domain containing the origin of the set of $z=u + i v \in \frac{1}{N}\Z^2$ such that $\{u' + i v' : |u-u'| \le (2N)^{-1}, |v-v'| \le (2N)^{-1} \}$ is contained in $D$. For each $N$, we get a measure on paths denoted $\mu^{\RW}_{\bd D_N}$ by considering the random walk excursion measure on $D_N$, and scaling the excursions by Brownian scaling: $\w^{(N)}(t) := N^{-1/2}\w(2Nt)$. As $N \to \infty$, these measures converge to $\mu_{\bd D}$, \defines{excursion measure on $D$}, which is an infinite measure on paths. Since Brownian motion in $\C$ is conformally invariant (up to a time-change), $\mu_{\bd D}$ is also conformally invariant. (See Proposition~\ref{confinvexcmeasD.nowprop}.) If $\Gamma$, $\Upsilon$ are disjoint arcs in $\bd D$, then conditioning the excursion measure to have endpoints $z\in \Gamma$, $w\in \Upsilon$, gives a probability measure on excursions from $\Gamma$ to $\Upsilon$ in $D$.

The primary result of this paper is that for any bounded, simply connected Jordan domain $D$, simple random walk excursion measure converges to Brownian excursion measure on $D$.

\begin{theorem}\label{final-thm}
If $D$ is a bounded, simply connected domain containing the origin with $\inrad(D)=1$, $\bd D$ is Jordan, and $D_N$ is the $1/N$-scale discrete approximation to $D$, then
$$\prohorov(\;4\,\mu^{\RW}_{\bd D_N}, \; \mu_{\bd D} \, ) \to 0$$
where $\prohorov$ denotes the Prohorov metric.
\end{theorem}

As we are discussing the convergence of infinite measures, we need to be a little careful about how we define convergence in the Prohorov metric (which is usually defined only for finite measures).  As the restriction of excursion measure to disjoint boundary arcs gives a finite measure, Theorem~\ref{final-thm} is to be interpreted as meaning that for \emph{any} pair of disjoint boundary arcs $\G$, $\U \subset D$,
$$\prohorov(4\mu^{\RW}_{\bd D_N}(\G_N, \U_N), \mu_{\bd D}(\G,\U)) \to 0$$
 where $\G_n$, $\U_N$ are the ``associated (discrete) boundary arcs in $D_N$.'' In Section~\ref{mainresultsection} we prove the precise formulation of Theorem~\ref{final-thm}.

Since a Brownian (resp., random walk) excursion can be viewed as consisting of a Brownian motion (resp., random walk) plus tails, the proof of convergence has two distinct parts---a ``global part'' plus a ``local part.''  The strong approximation of Koml\'os, Major, and Tusn\'ady~\cite{KomMT1, KomMT2} is used to couple random walk and Brownian motion in the interior of the domain away from the boundary. This global part does not depend on the smoothness of the boundary. The local part concerns the tails whose behaviour can be controlled using the Beurling estimates; here the structure of the boundary does come into play. The proof of convergence also employs an estimation of the discrete excursion Poisson kernel in terms of the excursion Poisson kernel derived in~\cite{KozL} which was used in that paper to prove a conjecture of Fomin~\cite{fomin}. Hence, by proving the weak convergence of excursion measures, we are extending the ``central limit theorem'' for the endpoints of the excursions proved in~\cite{KozL}. 

Technically, since $\mu_{\bd D_N}^{\RW}$ is supported on continuous curves, we must associate to $D_N$ a domain in $\C$ by identifying each point in $D_N$ with the square of side length $1/N$ centred at that point. It is important that these so-called ``union of squares'' domains $\tilde D_N$ converge to the original domain $D$. However, the convergence is not in the usual topological sense, but rather in the Carath\'eodory sense. This is captured by the following theorem which is carefully stated and proved in Section~\ref{cara-conv-section}. 

\begin{theorem}\label{cara-thm-intro}
If $f_N$, $f$ are conformal transformations of the unit disk $\D$ onto  $\tilde{D}_N$, $D$, respectively, with $f_N(0)=f(0)=0$ and $f'_N(0)$, $f'(0)>0$, then $f_N \to f$ uniformly on compact subsets of $\D$. In other words, $\tilde{D}_N \cara D$.
\end{theorem}

\subsection{Outline of the paper}

In Section~\ref{Sect2}, we establish some notation, and recall some facts from complex analysis about conformal transformations. We also review the definitions and basic facts about Green's functions on both $\C$ and $\Z^2$. Section~\ref{Sect3} is devoted to a discussion of excursions and excursion measures. Included are some fundamental ideas about spaces of curves and measures on metric spaces. We also review the Prohorov topology, and prove several easy lemmas about the Prohorov metric which are needed in the sequel. The Poisson kernel and excursion Poisson kernel are then reviewed, with an emphasis on their conformal covariance properties, and a construction of excursion measure on $D$, differing from that in~\cite{LawlerSLE}, is carried out. The final section, Section~\ref{Sect4}, is devoted to the proofs of Theorem~\ref{final-thm} and Theorem~\ref{cara-thm-intro}. The material in Section~\ref{scale-sect} relies on the results obtain in~\cite{KozL}. Instead of simply recopying those results as originally proved, we have translated them to statements in terms of  $D_N$, the $1/N$-scale discrete approximation to $D$. A review of some recent strong approximation results is included in Section~\ref{KMTsection} because of their necessity in the proof of Theorem~\ref{final-thm}. 

\section{Background and notation}\label{Sect2}

\subsection{Simply connected subsets of $\C$ and $\Z^2$}\label{scs}

A function $f:D \to D'$ is a \defines{conformal transformation} 
if $f$ is an analytic, univalent (i.e, one-to-one) function that is onto $D'$. 
It follows that $f'(z) \neq 0$ for $z\in D$, and $f^{-1}:D' \to D$
is also a conformal transformation. If $D$ is a domain in $\C$, then a connected $\Gamma \subset \bd D$ is an {\em (open) analytic arc} of $\p D$ if there
is a domain $E \subset \C$ that is symmetric
about the real axis  and a
conformal transformation $f: E \rightarrow f(E)$ such that
$f(E \cap \R) = \Gamma$ and $f(E \cap \Half) = f(E) \cap D$. We
say that $\p D$ is {\em locally analytic} at $x \in \p D$ if there
exists an analytic arc of $\p D$ containing $x$. For $D \subset \C$ with $0 \in D$,    
define the \defines{radius} (with respect to
 the origin) of $D$ to be $\rad(D)  := \sup\{|z| : z \in \bd D\}$,
 and the \defines{inradius} (with respect to the origin) of
 $D$ to be $\inrad(D) := \dist(0,\bd D) := \inf\{|z| : z \in \bd D\}$. 
The \defines{diameter} of $D$ is $\diam(D):=\sup\{|x-y| : x,y\in D\}$.
 Call a bounded domain $D$ a {\em Jordan domain}  if $\bd D$
 is a Jordan curve (i.e., homeomorphic to a circle).
  A Jordan domain is {\em nice} if the Jordan
 curve $\bd D$ can be expressed as a finite union of analytic curves. 
Note that Jordan domains  are necessarily simply connected.
For each $r>0$, let $\Dr{r}$ be the set of nice Jordan domains
 containing the origin of inradius $r$,
 and write $\scd := \bigcup_{r>0}\Dr{r}$. We also define 
$\scd^*$ to be the set of Jordan domains
containing the origin, and note that $\scd \subsetneq \scd^*$. If $D$, $D' \in \scd^*$,
 let $\mathcal{T}(D,D')$ be the set of all conformal transformations of $D$ onto $D'$. The Riemann mapping theorem implies that $\mathcal{T}(D,D') \neq \emptyset$, and since $\bd D$, $\bd D'$
are Jordan, the Carath\'eodory extension theorem tells us that $f \in \mathcal{T}(D,D')$ can be extended to a homeomorphism of $\overline{D}$ onto $\overline{D'}$. 
The statements and details of these two theorems may be found in~\cite[\S 1.5]{Duren}. 

There are three standard ways to define the boundary of a proper subset $A$ of $\Z^2$:
\begin{center}
\begin{itemize}
\item \textit{(outer) boundary:} $\;\bd A := \{y \in \Z^2 \setminus A: 
|y-x| = 1 \text{ for some } x \in A \}$;\\
\item \textit{inner boundary:} $\;\bd_i A :=\bd(\Z^2 \setminus A) =
 \{x\in A: |y-x|=1 \text{ for some } y \in \Z^2 \setminus A\}$;\\
\item \textit{edge boundary:} $\;\bd_e A := \{(x,y): x \in A,\, 
y \in \Z^2 \setminus A,\, |x-y| = 1\}$.
\end{itemize}
\end{center}

To each finite, connected $A \subset \Z^2$ we associate a domain  
$\tilde A \subset \C$ in the following way. For each edge $(x,y) \in \bd_e A$, considered as a 
line segment of length one, let $\ell_{x,y}$ be the perpendicular 
line segment of length one intersecting $(x,y)$ in the midpoint. 
 Let $\bd \tilde A$ denote the union of the line segments $\ell_{x,y}$,
 and let $\tilde A$ denote the domain with boundary $\bd \tilde A$ 
containing $A$. Observe that
\begin{equation*}
\tilde A \cup \bd \tilde A = \bigcup_{x \in A}
 \Square_x \;\text{ where }\; \Square_x := x + 
\left(\,[-1/2,1/2] \times [-1/2,1/2]\,\right).
\end{equation*}
That is, $\Square_x$ is the closed square of side 
length one centred at $x$ whose sides are parallel 
to the coordinate axes.  Also, note that $\tilde A$ is 
a simply connected domain if and only if $A$ is a simply
connected subset of $\Z^2$. We say $\tilde{A}$ is the
 \defines{``union of squares''} domain associated to $A$.

Let $\A$ denote the set of all finite simply connected 
subsets of $\Z^2$ containing the origin. If $A \in \A$, let 
$\inrad(A):=\min\{|z| : z \in \Z^2 \setminus A \}$ 
and $\rad(A):=\max\{|z| : z \in A\}$
denote the \defines{inradius} and \defines{radius} 
(with respect to the origin), respectively, of $A$, and 
define $\A^n$ to be the set of $A \in \A$ with $n \leq \inrad(A) \leq 2n$; 
thus $\A := \bigcup_{n> 0}\A^n$.  If $A \in \A$, $0 \neq x_1 \in \bd_i A$,
 and $[x_1,x_2,\ldots,x_j]$ is a nearest neighbour path in $A \setminus \{0\}$, then
the connected component of  $A \setminus  \{x_1,\ldots,x_j\}$
 containing the origin is simply connected. 

Finally, if $A \in \A$ with associated domain $\tilde A \subset \C$, 
then we write
$f_A := f_{\tilde{A}}$ for the conformal transformation of $\tilde A$ onto 
the unit disk $\D$ with $f_A(0) =0$, $f_A'(0) > 0$.

\subsection{Green's functions on $\C$ and $\Z^2$}\label{g-sect}

Let $D$ be a domain whose boundary contains a curve, and write
 $g_D(x,y)$ for the \defines{Green's function for Brownian
motion on $D$}. If $x \in D$, we can define $g_D(x, \cdot)$ 
as the unique harmonic function on $D\setminus \{x\}$, 
vanishing on $\bd D$, with $g_D(x,y) = -\log|x-y| + O(1)$ as  $|x-y| \to 0$. 
Equivalently, if $D \in \scd^*$, then for distinct points $x$, $y \in D$, 
$g_D(x,y) = \E^x[\log|B_{T_D}-y|] - \log |x-y|$
where $T_D := \inf \{t : B_t \not\in D \}$. 
In particular, if $0 \in D$, then $g_D(x) = \E^x[\log|B_{T_D}|] - \log |x|$  for  $x\in D$ 
where $g_D(x) :=g_D(0,x)$. For further details, consult~\cite[Chapter 2]{LawlerSLE}\label{SLE1}.
Since the Green's function is a well-known example of a conformal invariant
(see, e.g.,~\cite[\S 1.8]{Duren}), 
in order to determine $g_D$ for arbitrary $D \in \scd^*$, it is 
enough to find $f_D \in \mathcal{T}(D,\D)$. 
Conversely, suppose that $D \subset \C$ is a simply connected
 domain containing the origin with Green's function $g_D$.
The unique conformal transformation of $D$ onto $\Disk$ with $f_D(0) = 0,
f_D'(0) > 0$ can be written as  
\begin{equation}\label{riemannmap}
f_D(x) = \exp\{-g_D(x) + i\theta_D(x)\}.
\end{equation}
Note that $-g_D + i \theta_D$ is analytic in 
$D \setminus \{0\}$. Suppose that $A\in \A$, and that $g_A(x,y) := g_{\tilde A}(x,y)$. As explained in~\cite{KozL}, the exact form of the Green's function gives
\begin{equation*}
g_A(x,y)=\log \left|\frac{\overline{f_A(y)}f_A(x)-1}{f_A(y)-f_A(x)}\right|.
\end{equation*}
If we write $\theta_A : =\theta_{\tilde{A}}$, then~(\ref{riemannmap}) implies
$f_A(x) = \exp\{-g_A(x) + i\theta_A(x)\}$.
 
Let $S_n$ be a simple random walk on $\Z^2$, and let
 $A \subsetneq\Z^2$. If $\tau_A := \min \{j \ge 0 : S_j \not\in A\}$,
 then  we let
\begin{equation*}
G_A(x,y) :=  \sum_{j=0}^{\infty} \P^x \{S_j=y, \tau_A > j \}
\end{equation*}
denote the \defines{Green's function for random walk on $A$}, and set $G_A(x) := G_A(x,0) = G_A(0,x)$. 
Write $a$ for the potential kernel for simple random walk defined by
\[a(x)  :=  \sum_{j=0}^\infty \left[\P^0\{S_j=0\} - \P^x \{ S_j=0\}\right].\]
It is known~\cite[Theorem 1.6.2]{LawlerGreen} that as $|x| \to \infty$,
\begin{equation}\label{1.2}
a(x) = \frac{2}{\pi}\,\log|x| + k_0 +o(|x|^{-3/2})
\end{equation}
where $k_0 := (2\varsigma + 3\ln 2)/\pi$ and $\varsigma$ is Euler's constant, and that 
$G_A(x)=\E^x[a(S_{\tau_A})]-a(x)$ for  $x\in A$.
The error in~(\ref{1.2}) will suffice for our purposes even though stronger results are known; see~\cite{FukU1}. We also recall a uniform estimate for $G_A(x)$, and a relationship between the Green's functions $G_A$ and $g_A$ which is proved in~\cite[Theorem~1.2]{KozL}.

\begin{theorem}  \label{greentheoremB}
 If $A \in \A^n$, then
$ G_A(0) = -\frac 2 \pi \, \log f_A'(0) + k_0 + O(n^{-1/3}
 \, \log n)$. 
Furthermore, if $x \neq 0$, then
\begin{equation*}
G_A(x) = \frac{2}{\pi} \, g_A(x) + k_x + O(n^{-1/3} \,  \log n). 
\end{equation*}
where
\begin{equation}\label{defn-kx}
   k_x :=  k_0 + \frac 2 \pi \, \log |x| - a(x). 
\end{equation}
\end{theorem}

We conclude by defining what it means for two boundary arcs to be separated.  Note that separation is always defined in terms of distance in the unit circle.

\begin{definition}\label{defnsep.new}
Suppose that $A \in \A$ and $D \in \scd^*$.  Let $\Gamma_1$, $\Upsilon_1 \subset \bd_i A$ with $\overline{\Gamma_1} \cap \overline{\Upsilon_1} = \emptyset$, let $\Gamma_2$, $\Upsilon_2 \subset \bd D$ with $\overline{\Gamma_2} \cap \overline{\Upsilon_2} = \emptyset$, and write $\theta_1=\theta_A$, $\theta_2=\theta_D$.  The \defines{separation of $\Gamma_j$ and $\Upsilon_j$}, $j=1,2$, written $\sep(\Gamma_j,\Upsilon_j)$, is defined to be
\begin{equation}\label{defn-sep}
\sep(\Gamma_j,\Upsilon_j) :=  \inf \{|\theta_j(x)-\theta_j(y)| : x\in \Gamma_j, y \in \Upsilon_j \},
\end{equation}
and the \defines{spread of $\Gamma_j$ and $\Upsilon_j$}, written $\spr(\Gamma_j,\Upsilon_j)$, is defined to be
\begin{equation}\label{defn-spr}
\spr(\Gamma_j,\Upsilon_j) :=  \sup \{|\theta_j(x)-\theta_j(y)| : x\in \Gamma_j, y \in \Upsilon_j \}.
\end{equation}
If $\Gamma_1$, $\Upsilon_1 \subset \bd A$ instead, then~(\ref{defn-sep}) and~(\ref{defn-spr}) hold with $\theta_A$ extended to $\bd A$ in the natural way.
\end{definition}

\section{Excursions and excursion measure}\label{Sect3}

Much of this material may be found in~\cite{KozL} and in the recent book~\cite{LawlerSLE}. We repeat the relevant material here without proof in order to standardize our notation, and to remind the reader of the most important facts. We do, however, prove a number of useful lemmas about the Prohorov metric in Section~\ref{measuresonMS}.

\subsection{Metric spaces of curves}\label{curvespace}

A \defines{curve} $\gamma : I\to \C$ shall always mean a continuous mapping of an interval $I \subseteq [0, \infty)$ into $\C$. Let $\curves$ denote the set of curves $\gamma : [0,t_{\gamma}] \to \C$ where $0 < t_{\gamma} < \infty$, and write $\gamma[0,t_{\gamma}] := \{z \in \C : \gamma(t)=z \, \text{ for some } \, 0 \le t \le t_{\gamma}\}$ and similarly for $\gamma(0,t_{\gamma})$. There are three natural metrics that we will consider on $\curves$. Following~\cite{LawW1}, define the metric 
\begin{equation*}
\metricnoparam(\gamma, \gamma') := \inf_{\phi}\left[ \sup_{0\le s \le t_{\gamma}}|\, \gamma(s)-\gamma'(\phi(s))\, |\right]
\end{equation*}
where the infimum is over all increasing homeomorphisms $\phi:[0,t_{\gamma}] \to [0,t_{\gamma'}]$. Call $\tilde{\gamma}$ a \defines{reparameterization} of $\gamma \in \curves$ with parameterization $\phi$ if $\phi:[0,t_{\gamma}] \to [0,t_{\tilde{\gamma}}]$ is an increasing homeomorphism such that $\gamma(t) = \tilde{\gamma}(\phi(t))$ for each $0 \le t \le t_{\gamma}$. If $\tilde{\gamma}$ is a reparameterization of $\gamma$ under $\phi$, then $\gamma$ is a reparameterization of $\tilde{\gamma}$ under $\phi^{-1}$, and we write $\gamma \param \tilde{\gamma}$.  Finally, let $\curves^{*}$ be the set of equivalence classes of curves $\gamma \in \curves$ under the relation $\param$, so that the metric $\metricnoparam$ identifies curves which are equal modulo time reparameterization. In fact, it can be shown that $(\curves^{*}, \metricnoparam)$ is a complete metric space~\cite[Lemma~2.1]{AizB1}.

In order to account for the time parameterization, however, we let 
\begin{equation*}
\metricparam(\gamma, \gamma') := \inf_{\phi}\left[ \sup_{0\le s \le t_{\gamma}} \left\{ \,|\,\gamma(s)-\gamma'(\phi(s))\,| + |s-\phi(s)| \, \right\} \right]
\end{equation*}
where again the infimum is over all increasing homeomorphisms $\phi:[0,t_{\gamma}] \to [0,t_{\gamma'}]$. The metric $\metricparam$ does not identify curves which are equal modulo time reparameterization. A convenient choice of parameterization is $\phi(s) = t_{\gamma'}s/t_{\gamma}$. Define
\begin{equation*}
\metric(\gamma, \gamma') := \sup_{0\le s \le 1} |\, \gamma(t_{\gamma}s) - \gamma'(t_{\gamma'}s) \,| + |t_{\gamma} - t_{\gamma'}|
\end{equation*}
and note that it is straightforward to verify $\metric$ is also a metric on $\curves$.  Neither $(\curves, \metricparam)$ nor $(\curves, \metric)$ is complete as Example~\ref{curvesnotcomplete} combined with the next lemma will show. For the proof of this lemma, consult~\cite[Lemma~5.1]{LawlerSLE}.

\begin{lemma}\label{metricequiv}
If $\gamma_1$, $\gamma_2 \in \curves$, then $\metricparam(\gamma_1, \gamma_2) \le \metric(\gamma_1, \gamma_2) \le \metricparam(\gamma_1, \gamma_2) + \osc(\gamma_2, 2\metricparam(\gamma_1, \gamma_2))$, where $\osc(\gamma,\delta) := \sup \{|\gamma(t) - \gamma(s)| : |t-s|\le \delta\}$ denotes the modulus of continuity of $\gamma$.
\end{lemma}

To account for the incompleteness of $(\curves, \metric)$, we consider a larger complete space $\BX$, and identify subspaces of $(\curves, \metric)$ with closed subspaces of $\BX$. Let $C[0,1]$ denote the space of continuous complex-valued functions on $[0,1]$ under the metric $d_{\infty}(\gamma^*_1, \gamma^*_2) := \sup_{0 \le r \le 1}|\gamma^*_1(r) - \gamma^*_2(r)|$, and denote the usual metric on $\R$ by $\abs$. Consider the separable Banach space $\BX := C[0,1] \times \R$ with metric $d_{\BX} := d_{\infty} + \abs$.  Thus, elements of $\BX$ are pairs $(\gamma^*, t)$ where $\gamma^* \in C[0,1]$, $t \in \R$, and $d_{\BX}((\gamma^*_1, s),(\gamma^*_2, t)) = \sup_{0 \le r \le 1}|\gamma^*_1(r) - \gamma^*_2(r)| + |s-t|$. We can embed $\curves$ into $\BX$ via $\iota : \curves \hookrightarrow \BX$, $\gamma \mapsto (\gamma^*, t_{\gamma})$, where $\gamma^*(r) := \gamma(t_{\gamma}r)$, $0 \le r \le 1$. However, $\iota(\curves) = \{(\gamma^*,t) \in \BX : t>0 \} =: \BX^+$ is not a closed subspace of $\BX$. The metric spaces $(\BX^+, d_{\BX})$ and $(\curves, d_{\BX,\curves})$ are isomorphic, where $d_{\BX, \curves}$ is the induced metric in $\curves$ associated to the metric $d_{\BX}$ in $\BX$. That is, if $\gamma_1$, $\gamma_2 \in \curves$, then $\iota(\gamma_i) = (\gamma_i^*, t_{\gamma_i})$, $i=1$, $2$, so that $d_{\BX,\curves}(\gamma_1, \gamma_2) = d_{\BX}((\gamma_1^*, t_{\gamma_1}), (\gamma_2^*, t_{\gamma_2}))$. It follows that $d_{\BX,\curves} = \metric$ and $(\curves, \metric) \cong (\BX^+, d_{\BX})$ since
\begin{equation*}
d_{\BX}((\gamma_1^*, t_{\gamma_1}), (\gamma_2^*, t_{\gamma_2})) = \sup_{0 \le r \le 1}|\gamma_1(t_{\gamma_1}r) - \gamma_2(t_{\gamma_2}r)| + |t_{\gamma_1}-t_{\gamma_2}| = \metric(\gamma_1, \gamma_2)
\end{equation*}

\begin{example}\label{curvesnotcomplete}
Suppose $\gamma \in \curves$ is given by $\gamma(r) = r + ir$, $0 \le r \le 1$, and for $n=1,2,\ldots$, let $\gamma_n(r) = nr + inr$, $0 \le r \le 1/n$. Notice that $\gamma_n^*=\gamma^*=\gamma$. Thus, $\iota(\gamma_n) = (\gamma_n^*,t_{\gamma_n}) = (\gamma^*, 1/n)$ so clearly $\{ (\gamma_n^*, t_{\gamma_n})\}$ is a Cauchy sequence in $\BX$, and  $\{\gamma_n\}$ is a Cauchy sequence in $(\curves,\metric)$. Since $\BX$ is complete, it has a limit, namely $(\gamma^*,0) \in \BX$.  However, $(\gamma^*,0) \not\in \BX^+=\iota(\curves)$ so that $(\gamma^*,0)$ does not have a counterpart in $\curves$. This shows that $(\curves, \metric)$ is not complete, and illustrates the reason for considering $\BX$.  
\end{example}

However, if the limit does have a counterpart in $\curves$ (i.e., if $(\gamma^*,t) \in \BX^+$ so that $\iota^{-1}(\gamma^*,t) \in \curves$), then we have the following result. See~\cite[Lemma~5.2]{LawlerSLE} for the proof.

\begin{lemma}\label{metrictopologyequiv.nowlem}
Let $(\gamma_n^*, t_n) \in \BX^+$ for $n=1,2, \ldots$, so that 
$\gamma_n:=\iota^{-1}(\gamma_n^*, t_n) \in \curves$. Suppose that for some $(\gamma^*,t) \in \BX$, 
$d_{\BX}((\gamma_n^*,t_n), (\gamma^*,t)) \to 0$. If $t>0$ so that $(\gamma^*,t) \in \BX^+$, then $\gamma :=\iota^{-1}(\gamma^*,t) \in \curves$, and $d_{\BX}((\gamma_n^*,t_n), (\gamma^*,t)) \to 0$ if and only if $\metricparam(\gamma_n, \gamma) \to 0$ as $n \to \infty$, or equivalently, $\metric(\gamma_n, \gamma) \to 0$ if and only if $\metricparam(\gamma_n, \gamma) \to 0$ as $n \to \infty$.
\end{lemma}

Consequently, $\metricparam$ and $\metric$ generate the same topology on $\curves$. Thus, when we need to discuss convergence or continuity in $\curves$, it can be with respect to whichever metric is more convenient for the given problem.

If $a >0$, let $\curves_a := \{\gamma \in \curves  : t_{\gamma} \ge a \}$, and set $\iota(\curves_a)=\{(\gamma^*,t) \in \BX : t \ge a\} =: \BX_a$.  Note that $\BX_a$ \emph{is} a closed subspace of $\BX$ so that $(\BX_a, d_{\BX}) \cong (\curves_a, \metric)$ is complete. However, $\curves_a$ is not complete under $\metricparam$. As an example, consider $\gamma_n(r) = r^n$, $0 \le r \le 1$, which is a Cauchy sequence in $(\curves_1, \metricparam)$ that has no limit.  By Lemma~\ref{metricequiv}, if $\{\gamma_n\}$ is a 
Cauchy sequence in $(\curves_a, \metricparam)$ that is equicontinuous, then it is a Cauchy sequence in $(\curves_a, \metric)$ and therefore has a limit. In what follows, we will refer to spaces of curves which are primarily subspaces of $\curves$.  Since such spaces are isomorphic to subspaces of $\BX$, we prefer to work with $(\curves,\metric)$ rather than $(\BX,d_{\BX})$ unless it is necessary to explicitly mention this isomorphism.

If $D$ is a simply connected proper subset of $\C$, and $\gamma \in \curves$, then we say that \defines{$\gamma$ is in $D$} if $\gamma(0,t_{\gamma}) \subset D$.  This does not require that either $\gamma(0) \in D$ or $\gamma(t_{\gamma}) \in D$. We define the space $\curves(D)$ as $\curves(D) := \{ \gamma \in \curves : \gamma \text{ is in } D \}$.
For $z$, $w \in \overline{D}$, let $\curves_z(D)$ be the set of $\gamma \in \curves(D)$ with $\gamma(0)=z$, let $\curves^w(D)$ be the set of $\gamma \in \curves(D)$ with $\gamma(t_{\gamma})=w$, and define $\curves_z^w(D) := \curves_z(D) \cap \curves^w(D)$.  Finally, if $\Gamma$, $\Upsilon \subset \bd D$ with $\overline{\Gamma} \cap \overline{\Upsilon} = \emptyset$, write 
$\curves_{\Gamma}^{\Upsilon}(D) := \bigcup_{z\in \Gamma, w\in\Upsilon}\curves_z^w(D)$.

\begin{definition}
Suppose $\gamma \in \curves(D)$. We say \defines{$\gamma$ is an excursion in $D$} if $\gamma(0) \in \bd D$ and $\gamma(t_{\gamma}) \in \bd D$, and
we say \defines{$\gamma$ is an excursion from $z$ to $w$ in $D$} if $\gamma(0)=z\in \bd D$ and $\gamma(t_{\gamma})=w  \in \bd D$, i.e., if $\gamma \in \curves_z^w(D)$ with $z$, $w \in \bd D$. If $\Gamma$, $\Upsilon \subset \bd D$ with $\overline{\Gamma} \cap \overline{\Upsilon} = \emptyset$, then we say \defines{$\gamma$ is a $(\Gamma, \Upsilon)$-excursion in $D$} if $\gamma(0)\in \Gamma$ and $\gamma(t_{\gamma}) \in \Upsilon$, i.e., if $\gamma \in \curves_{\Gamma}^{\Upsilon}(D)$.
\end{definition}

Suppose that both $D$ and $D'$ are simply connected domains in $\C$, and  $f:D \to D'$ is a conformal transformation. For $\gamma \in \curves(D)$, let
$$A_s=A_{s,f,\gamma} := \int_0^s |f'(\gamma(r))|^2 \, \d r \;\;\text{ and }\;\; \sigma_t  = \sigma_{t,f,\gamma} := \inf\{s : A_s \ge t\}.$$ 
If $\gamma \in \curves(D)$ with $A_{t_{\gamma}} < \infty$, and if $f$ extends to the endpoints of $\gamma$, then we define the \defines{image\label{imagepage} of $\gamma$ under $f$}, denoted $f \circ \gamma \in \curves(D')$, by setting $t_{f \circ \gamma} := A_{t_{\gamma}}$ and $f\circ \gamma (t) := f(\gamma(\sigma_t))$ for $0 \le t \le A_{t_{\gamma}}$. Since $s \mapsto A_{s, f,{\gamma}}$ is non-negative, continuous, and strictly increasing, it follows that $t \mapsto \sigma_{t, f,{\gamma}}$ is well-defined. The following is a special case.

\begin{example}\label{brscale_example}
Let $D$ be a simply connected proper subsets of $\C$, and for $a \in \C \setminus \{0\}$, let $f_a(z)=az$.  If $\gamma \in \curves(D)$, then we define the \defines{Brownian scaling map} $\brscale_a : \curves(D) \to \curves(f_a(D))$ by setting $t_{\brscale_a \gamma}:=|a|^2t_{\gamma}$ and $\brscale_a \gamma(t) := a \, \gamma\left(|a|^{-2}t\right)$ for $0 \le t \le t_{\brscale_a \gamma}$.
\end{example}

In  particular, if $D$, $D' \in \scd$, $\gamma$ is an excursion in $D$, and $f \in \mathcal{T}(D,D')$ so that $f$ \emph{does} extend to the endpoints of $\gamma$, then $f \circ \gamma =: \gamma' \in \curves(D')$ is an excursion in $D'$. Note that $t_{\gamma'} = A_{t_\gamma}$ (i.e., $\sigma_{t_{\gamma'}}=t_{\gamma}$) and $\gamma'(t) = f(\gamma(\sigma_t))$ for $0 \le t \le t_{\gamma'}$.

\begin{example}
As an application of Brownian scaling, suppose that $f(z)=(1+\eps)z$ for $z \in \D$, $0 < \e <1$, and let $\gamma$ be an excursion from $x$ to $y$ in $\D$.  Then $\gamma' := f \circ \gamma$ is an excursion from $(1+\eps)x$ to $(1+\eps)y$ in $(1+\eps)\D$ given explicitly by $\gamma'(t) = (1+\eps)\,\gamma\left((1+\eps)^{-2}t\right)$ for $0 \le t \le t_{\gamma'}=(1+\eps)^2t_{\gamma}$. Furthermore, it is not very difficult to verify that there exists a constant $C=C(\gamma)$ such that $\metric(\gamma, \gamma') \le C \eps$.
\end{example}

If $E$ is a domain containing $\overline{D}$ and $f$ is a conformal mapping of $E$, then it follows from the Koebe growth and distortion theorems~\cite[Theorems 2.4, 2.5, 2.6]{Duren} that $|f'|$, $|f''|$, and $1/|f'|$ are uniformly bounded on $D$, and the function $\gamma \mapsto f \circ \gamma$ from $\curves(D)$ to $\curves(f(D))$ is continuous.  If $D \in \scd$, then since $\bd D$ is piecewise analytic, $\bd D = \bigcup_{i=1}^n\Gamma_i$ for some finite union of analytic curves $\G_i$. Hence, any conformal mapping $f$ of $D$ can be analytically continued across each $\G_i$ so 
$\gamma \mapsto f \circ \gamma : \curves(D) \to \curves(f(D))$ is continuous; we denote this induced map by $f$.

\begin{definition}
If $\gamma_1$, $\gamma_2 \in \curves$ with $\gamma_1(t_{\gamma_1}) = \gamma_2(0)$, then we define the \defines{concatenation of $\gamma_1$ and $\gamma_2$}, denoted $\gamma_1 \oplus \gamma_2$, by setting $t_{\gamma_1 \oplus \gamma_2} := t_{\gamma_1} + t_{\gamma_2}$, and
\begin{equation*}
\gamma_1 \oplus \gamma_2 (t) := 
\begin{cases}
\gamma_1(t),                    &\text{if } 0 \le t \le t_{\gamma_1},\\
\gamma_2(t-t_{\gamma_1}),       &\text{if } t_{\gamma_1} \le t \le t_{\gamma_1 \oplus \gamma_2}. 
\end{cases}
\end{equation*}
\end{definition}

Note that $(\gamma_1, \gamma_2) \mapsto \gamma_1 \oplus \gamma_2$ is a continuous map from $\curves^{w} \times \curves_w$ to $\curves$ for every $w \in \C$.

\begin{definition}
If $0\le r < s \le t_{\gamma}$, then we define the \defines{truncation operator} $\truncopr{r}{s} : \curves \to \curves$ by setting $t_{\trunc{r}{s}} := s-r$ and $\trunc{r}{s}(t)  :=\gamma(r+t)$ for $0 \le t \le t_{\trunc{r}{s}}$.
\end{definition}

Observe that $\trunc{r}{s}[0, t_{\trunc{r}{s}}] = \gamma[r,s]$, and that by definition, truncation undoes concatenation. If $\gamma_1$, $\gamma_2 \in \curves$ with $\gamma_1(t_{\gamma_1}) = \gamma_2(0)$, then
$\truncopr{0}{t_{\gamma_1}} \gamma_1 \oplus \gamma_2(t)= \gamma_1(t)$, $0 \le t \le t_{\gamma_1}$, and 
$\truncopr{t_{\gamma_1}}{t_{\gamma_1\oplus \gamma_2}} \gamma_1 \oplus \gamma_2(t)= \gamma_2(t)$, $0 \le t \le t_{\gamma_2}$. It is easily seen that 
$$\metricparam(\trunc{r}{s}, \gamma) \le r + (t_{\gamma}-s) + \diam(\gamma[0, r]) + \diam(\gamma[s, t_{\gamma}]).$$
Therefore, if $r_n \to 0+$ and $s_n \to t_{\gamma}-$, then by Lemma~\ref{metrictopologyequiv.nowlem}, $\metric(\trunc{r_n}{s_n}, \gamma)\to 0$.

\subsection{General facts about measures on metric spaces}\label{measuresonMS}

Throughout this section, suppose that $(\Xi, \rho)$ is a metric space. Let $\mathcal{B}_{\rho} := \mathcal{B}_{\rho}(\Xi)$ denote the Borel $\sigma$-algebra associated to the topology induced by $\rho$, so that $(\Xi, \mathcal{B}_{\rho})$ is a measurable space. A measure $m$ on $(\Xi, \rho)$ will always be a $\sigma$-finite measure on $(\Xi, \mathcal{B}_{\rho})$. Denote the \defines{total mass} of $m$ by $|m| := m(\X)$.  If $|m| < \infty$, then  $m$ is a \defines{finite measure}; otherwise it is an \defines{infinite measure}. Denote the space of finite (resp., probability) measures on $(\Xi, \mathcal{B}_{\rho})$ by $\M(\Xi)$ (resp., $\probmeas(\X)$). If $m \in \finmeas(\X)$ with $|m| >0$, we write $m^{\#} : = m/|m|$ so that $m^{\#} \in \probmeas(\X)$. Recall (see~\cite{Bil1}) that every finite measure $m$ on $(\X, \mathcal{B}_{\rho})$ is \defines{regular}; i.e., if $V \in \mathcal{B}_{\rho}$ and $\eps > 0$, then there exist a closed set $F$ and an open set $G$ such that $F \subseteq V \subseteq G$ and $m(G\setminus F) < \eps$. If $V \in \mathcal{B}_{\rho}$, then we say that $m$ is \defines{concentrated on $V$} if $m(\X \setminus V ) = 0$. Observe that $V$ need not be closed.  

\begin{definition}\label{defn-prohorov}
If $m_1$, $m_2 \in \M=\M(\Xi)$, let $\prohorov : \M \times \M \to [0,\infty)$ denote the \defines{Prohorov metric} given by $\prohorov(m_1,m_2) := \inf\{\eps>0 : m_1(F) \le m_2(F^{(\eps)}) + \eps, \;\; m_2(F) \le m_1(F^{(\eps)}) + \eps \;\; \forall \; F \in \mathcal{B}_{\rho} \}$ where $F^{(\eps)} := \{x\in \Xi : \rho(x,y) < \eps, \text{ some } y\in F\}$. 
\end{definition}

It is easily checked that $(\finmeas(\Xi), \prohorov)$ is itself a metric space. Observe that $F^{(\eps)}$ is Borel, and that symmetry follows since $((F^{(\eps)})^c)^{(\eps)} \subseteq F^c$. If $m_1$, $m_2 \in \probmeas(\Xi)$, then an equivalent definition of $\prohorov$ is given by
$\prohorov(m_1,m_2) = \inf\{\eps>0 : m_1(F) \le m_2(F^{(\eps)}) + \eps \text{ for every closed } F \in \mathcal{B}_{\rho} \}$. It is known~\cite[Theorem~2.4.2]{Bor1} that both metrics on $\probmeas(\Xi)$ are equivalent and consistent with the Prohorov topology. Also note that $|\,|m_1|-|m_2|\,| \le \prohorov(m_1, m_2) \le \max \{ |m_1|, |m_2|\}$.  The following two theorems are standard.

\begin{theorem}\label{bergstromthm.nowprop.nowthm}
Suppose that $f$ is a continuous mapping of the metric space $(\X, \rho)$ into the metric space $(\X', \rho')$.  Then a measure $m$ on $(\X, \mathcal{B}_{\rho})$ determines a measure $m'$ on $(\X', \mathcal{B}_{\rho'})$ such that $f\circ m (V') = m'(V') = m(f^{-1}(V'))$ for any $V' \in \mathcal{B}_{\rho'}$. That is, $f\circ m \in \finmeas(\X')$ is given explicitly by $f \circ m (V') := m (\{x\in \X : f(x) \in V'\})$ for any $V' \in \mathcal{B}_{\rho'}$. Furthermore,
$$\int_{\Xi'}h(x') \, m'(\d x') = \int_{\Xi}h(f(x)) \, m(\d x)$$
for any bounded, continuous function $h: \Xi' \to \R$.  
\end{theorem}

\begin{theorem}\label{weak_is_prohorov}
If $(\Xi, \rho)$ is a complete, separable metric space, then the metric space $(\probmeas(\X), \prohorov)$ is also complete and separable, where $\prohorov$ is the Prohorov metric as in Definition~\ref{defn-prohorov}. Furthermore, if $m_n$, $m \in \probmeas(\X)$, then as $n \to \infty$, $\prohorov(m_n, m) \to 0$ if and only if $m_n \Rightarrow m$ weakly. 
\end{theorem}

\begin{importantremark}\label{important-remark}
Whenever we say a sequence of measures converges, it will be with respect to the Prohorov metric.
\end{importantremark}

As noted in~\cite[page 29]{Bor1}, Strassen proved another equivalent definition of $\prohorov$ consistent with the Prohorov topology is given by
\begin{equation*}
\prohorov(m_1, m_2) = \inf_{\mathfrak{M}} \,\left[\,\inf \{ \eps>0 : \P\{\rho(X_1,X_2) \ge\eps \} \le \eps \} \,\right],
\end{equation*}
where $\mathfrak{M}$ is the set of all $\Xi \times \Xi$-valued random variables $(X_1,X_2)$ with $\mathcal{L}(X_1)=m_1$ and  $\mathcal{L}(X_2)=m_2$ where $\mathcal{L}$ denotes law. In fact, an easy calculation shows that if $X_i$ are $(\Xi, \rho)$-valued random variables with $\mathcal{L}(X_i) =m_i$, $i=1, 2$, and if $\P\{\rho(X_1,X_2) \ge\eps \} \le \eps$, then $\prohorov(m_1,m_2) \le \eps$. If $(\X, \rho)$ is a complete and separable metric space, then to show a sequence of non-zero finite measures $m_n \in \finmeas(\X)$ converges to $m \in \finmeas(\X)$, it suffices to show that both $|m_n| \to |m|$ and $\prohorov(m_n^{\#}, m^{\#}) \to 0$, as $n \to \infty$. In particular, we record the following version of these remarks.

\begin{proposition}\label{d-p-prop}
Let $\gamma$, $\gamma'$ be $\curves$-valued random variables with $\mathcal{L}(\gamma) = \mu$ and $\mathcal{L}(\gamma') = \mu'$, respectively.  If $\P\{\metric(\gamma,\gamma') \ge \eps\} \le \eps$, then $\prohorov(\mu,\mu')\le \eps$.
\end{proposition}

\begin{lemma}\label{thankjose1.nowlem}
Suppose that $(\X,\rho)$ is a complete, separable metric space, and that $m_1$, $m_2 \in \finmeas(\X)$.  If $C >0$, then 
$\prohorov(Cm_1,Cm_2) \le (C\vee 1)\, \prohorov(m_1,m_2)$.
\end{lemma}

\begin{proof}
Suppose that $\prohorov(m_1,m_2) = \eps$. To begin, let $C>1$. Then since $m_1(F) \le m_2(F^{(\eps)}) + \eps$ for every $F$ Borel, we have $Cm_1(F) \le Cm_2(F^{(\eps)}) + C\eps$.  Since $C\eps > \eps$, we have $F^{(\eps)} \subset F^{(C\eps)}$. Hence, $Cm_1(F) \le Cm_2(F^{(C\eps)}) + C\eps$. Interchanging $m_1$ and $m_2$ gives $\prohorov(Cm_1,Cm_2)  \le C\eps$. Suppose instead that $C<1$. Then since $m_1(F) \le m_2(F^{(\eps)}) + \eps$, and $C\eps < \eps$, we have  $m_1(F) \le m_2(F^{(\eps)}) + \eps/C$.  Multiplying by $C$ gives $Cm_1(F) \le Cm_2(F^{(\eps)}) + \eps$. Interchanging $m_1$ and $m_2$ yields $\prohorov(Cm_1,Cm_2)  \le \eps$. Thus, the conclusion follows.
\end{proof}

\begin{lemma}\label{fofcm.prop.nowlem}
If $f: (\X,\rho) \to (\X',\rho')$ is continuous, $m \in \finmeas(\X)$, and $C$ is a constant, then 
\begin{equation}\label{fofcm.eq}
f\circ(Cm) = C(f\circ m).
\end{equation}
\end{lemma}

\begin{proof}
Since $f\circ(Cm)(V') = (Cm)(f^{-1}(V')) = C\left[m(f^{-1}(V'))\right] = C\left[f\circ m (V')\right]$ for any $V' \in \mathcal{B}_{\rho'}(\X')$ the  result follows. 
\end{proof}

\begin{lemma}\label{feb8.lem1.nowprop.nowlem}
Suppose that $(\X, \rho)$ is a complete, separable metric space, and let $m \in \M(\X)$.  If $f_n, f: (\X, \rho) \to (\X, \rho)$ are continuous, and $f_n \to f$ uniformly, then 
$\prohorov(f_n \circ m, f\circ m) \to 0$. 
\end{lemma}

\begin{proof}
Assume first that $m \in \probmeas(\X)$. If $\mu_n := f_n \circ m$ and $\mu:=f\circ m$, then by Theorem~\ref{weak_is_prohorov}, it suffices to show that $\mu_n \Rightarrow \mu$ weakly. Suppose that $h:\X \to \R$ is a bounded, continuous function.
Hence, by Theorem~\ref{bergstromthm.nowprop.nowthm}, we conclude that 
$$\int_{\X} h(x) \,\mu_n(\d x) = \int_{\X} h\circ f_n (x) \,m(\d x) \to \int_{\X} h \circ f (x) \,m(\d x) = \int_{\X} h(x) \,\mu(\d x)$$
since $f_n \to f$ uniformly. We next consider $m \in \finmeas(\X)$. If $m$ is the zero measure, the result is trivial. If $|m|>0$, then by~(\ref{fofcm.eq}) and Lemma~\ref{thankjose1.nowlem}, $\prohorov(f_n \circ m, f\circ m) = \prohorov(|m|\,(f_n \circ m^{\#}), |m|\,(f\circ m^{\#})) =(|m| \vee 1)\, \prohorov(f_n \circ m^{\#}, f\circ m^{\#}) \to 0$. 
\end{proof}

\begin{lemma}\label{thankjose3.nowcor.nowlem}
Under the same assumptions as Lemma~\ref{feb8.lem1.nowprop.nowlem}, if $m_2 \in \M(\X)$, and $\prohorov(f_n \circ m_1, f\circ m_1) \to 0$ as $n \to \infty$, then $\prohorov(f_n \circ m_1, m_2) \to \prohorov(f \circ m_1, m_2)$.
\end{lemma}

\begin{proof}
Since $\prohorov(\cdot, \cdot)$ is a metric, we have by the triangle inequality
$\prohorov(f_n \circ m_1, m_2) \le \prohorov(f_n \circ m_1, f\circ m_1) + \prohorov(f\circ m_1, m_2)$,
so that
\begin{equation}\label{cor.eq1}
\limsup_{n \to \infty} \prohorov(f_n \circ m_1, m_2) \le \prohorov(f\circ m_1, m_2).
\end{equation}
However, the reverse inequality tells us that 
$\prohorov(f \circ m_1, m_2) \le \prohorov(f \circ m_1, f_n\circ m_1) + \prohorov(f_n \circ m_1, m_2)$, 
so that 
\begin{equation}\label{cor.eq2}
\liminf_{n \to \infty} \prohorov(f_n \circ m_1, m_2) \ge \prohorov(f\circ m_1, m_2).
\end{equation}
By combining~(\ref{cor.eq1}) and~(\ref{cor.eq2}), the result follows.
\end{proof}

We conclude this section by reviewing how to define a measure by Riemann integration. Let $\Lambda \subset \C$ be an analytic arc that is parameterized by $\xi : [0,t_{\xi}] \to \C$ with $t_{\xi}<\infty$. Consider the measures $\{ \mu(z, \cdot) : z \in \Lambda\} \subset \finmeas(\X)$, and let $\{\xi(0)=z_0, z_1, \ldots, z_n=\xi(t_{\xi})\}$ partition $\Lambda$. Let $z_i^* \in [z_{i-1}, z_i]$, $|\Delta z_i| = |z_{i}-z_{i-1}|$, $i=1, \ldots n$,  and set
$$\mu_n(\cdot) := \sum_{i=1}^n \mu(z_i^*, \cdot) \, |\Delta z_i|.$$
Note that $\mu_n(\cdot) \in \finmeas(\X)$ for each $n$. If $\lim_{n\to \infty} \mu_n(\cdot)$ exists in $\finmeas(\X)$, then we define the Riemann integral of the measure-valued function $z \in \Lambda \mapsto \mu(z, \cdot) \in \finmeas(\X)$ to be this limiting value; that is,
\begin{equation}\label{Rlimit}
\mu(\cdot) := \int_{\Lambda} \mu(z, \cdot) \, |\d z| := \lim_{n\to \infty} \mu_n(\cdot).
\end{equation}
Several conditions guarantee the existence of the Riemann integral. For instance,  if $z \mapsto \mu(z, \cdot)$ is continuous at $z_0$ for all $z_0 \in \Lambda$, or if $(\Xi, \rho)$ is a complete and separable metric space and $\{ \mu_n(\cdot)\}$ is a Cauchy sequence,  then~(\ref{Rlimit}) exists in $\finmeas(\Xi)$.

\subsection{Excursion Poisson kernel}\label{EPKnewsection}

We now briefly review several results about the Poisson kernel and the excursion Poisson kernel. Further details including proofs can be found in~\cite{KozL}. Let $D$ be a domain, and write $\mathbf{n}_x = \mathbf{n}_{x,D}$ for the unit normal at $x$ pointing into $D$. If $x \in D$ and $\bd D$ is locally analytic at  $y \in \bd D$, then both \defines{harmonic measure} and its density with 
respect to arc length, the \defines{Poisson kernel} $H_D(x,y)$, are well-defined. The behaviour of the Poisson kernel under a conformal transformation can be easily deduced from the Riemann mapping theorem and L\'evy's theorem on the conformal invariance of Brownian motion~\cite{Bass}. See also~\cite[Proposition~2.10]{KozL}.

\begin{proposition}\label{PKconf}
If $f: D \to D'$ is a conformal transformation,
$x \in D$,  $\p D$ is locally analytic at $y \in \bd D$,
and $\p D'$ is locally analytic at $f(y)$,
then
\begin{equation}\label{conf-inv-PK}
H_D(x,y) = |f'(y)| \, H_{D'}(f(x), f(y)).
\end{equation}
Furthermore, if $\Gamma \subset \bd D$ and $f \circ \Gamma \subset \bd D'$ are analytic, then
\begin{equation}\label{PKjune05}
H_D(x,\Gamma) := \int_{\Gamma} H_D(x,y) \; |\d y| = H_{D'}(f(x), f\circ \G).
\end{equation}
\end{proposition}

\begin{definition}
For $x$, $y \in \bd D$, $x \neq y$, the \defines{excursion Poisson kernel} 
$H_{\bd D}(x,y)$ is given by
\[    H_{\bd D}(x,y) := \lim_{\epsilon \rightarrow 0+}
                  \frac{1}{\epsilon} \, H_{  D}(x + \eps\mathbf{n}_x,y)
                   = \lim_{\epsilon \rightarrow 0+}
                  \frac{1}{\epsilon} \, H_{  D}( y + \eps\mathbf{n}_y,x). \]
\end{definition}

For a proof of the next proposition, see~\cite[Proposition~2.11]{KozL}.

\begin{proposition}\label{EPKconf}
Suppose   $f: D \to D'$ is a conformal 
transformation and $x$, $y$ are distinct points
on $\p D$.  Suppose that $\bd D$ is locally analytic
at $x$, $y$, and $\bd D'$ is locally analytic
at $f(x)$, $f(y)$.
Then
$H_{\bd D}(x,y) = |f'(x)| \,|f'(y)| \, H_{\bd D'}(f(x), f(y))$.
\end{proposition}

\begin{corollary}\label{oneprop.nowcor}
If $x$, $y \in \bd \D$, $x \neq y$, and $f \in \mathcal{T}(\D,\D)$ with $f(x)=x$ and $f(y)=y$, then $|f'(x)|\,|f'(y)|=1$.
\end{corollary}

\begin{proof}
If $f(x)=x$ and $f(y)=y$, then we immediately obtain from Proposition~\ref{EPKconf} that
$H_{\bd \D}(x,y) = |f'(x)| \, |f'(y)| \, H_{\bd \D}(f(x),f(y)) = |f'(x)| \, |f'(y)| \, H_{\bd \D}(x,y)$. 
\end{proof}

\subsection{Brownian excursion measures on $(\curves, \metric)$}

We now remind the reader of several Brownian measures on $(\curves, \metric)$, and outline the construction of the Brownian excursion measure on $D$. The exposition follows~\cite{LawlerSLE}, although there are some noticeable differences. We begin with a general definition.

\begin{definition}
A \defines{measure $\mu$ on $\curves$} is defined to be a $\sigma$-finite measure on the measurable space $(\BX, \B_{d_{\BX}})$ concentrated on $\BX^+$. 
\end{definition}

Suppose $B_t$ is a Brownian motion with $B_0=z$, and let $T_D := \inf \{t : B_t \not\in D\}$ be its exit time from $D$.  The process $X_t^D := B_{t \wedge T_D}$, $t \ge 0$, is Brownian motion killed on exiting $D$. Let $D \in \scd$ and suppose $w \in \bd D$ so that the Poisson kernel $H_D(z,w)$ is well-defined. Define the continuous, positive martingale $M$ by $M_t^D := H_D(X_t^D,w) / H_D(z,w)$, and the probability  $\P^{z,w}_D$ by $\P^{z,w}_D(A) := \E^z[M_t^D; A]$ for $A \in \F_t$. As noted in~\cite{Bass}, the law of the process $X_t^D$ under $\P^{z,w}_D$ is that of Brownian motion conditioned to exit $D$ at $w$. 

\begin{definition}\label{int-bnd-exc-meas}
Suppose that $D \in \scd$.  The \defines{interior-to-boundary excursion measure from $z$ to $w$ in $D$}, written $\mu_D(z,w)$, is defined to be $\mu_D(z,w) := H_D(z,w) \cdot \P^{z,w}_D$,
and the \defines{interior-to-boundary excursion measure from $z$ in $D$}, written $\mu_D(z)$, is defined by
\begin{equation}\label{wienerexc}
\mu_D(z) := \int_{\bd D} \mu_D(z,w) \, |\d w|.
\end{equation}
\end{definition}

Translation invariance and Brownian scaling imply $w \mapsto \mu_D(z,w)$ is continuous so that~(\ref{wienerexc}) is well-defined as in~(\ref{Rlimit}). The measure on paths $\mu_D(z)$ is what is generally called \defines{Wiener measure}. Observe that $\mu_D(z)$ is a measure on $\curves$ concentrated on $\curves_z(D)$, and that by definition, $\mu_D(z,w)$ is a finite measure with mass $|\mu_D(z,w)| = H_D(z,w)$. As such we can consider the normalized probability measure 
\begin{equation}\label{normalized-int-bnd-exc-meas}
\mu^{\#}_{D}(z,w) := \frac{\mu_D(z,w)}{|\mu_D(z,w)|}= \frac{\mu_D(z,w)}{H_D(z,w)} := \P^{z,w}_D.
\end{equation}

It is well-known that two-dimensional Brownian motion is conformally invariant, and consequently so too is Wiener measure. 
We express this as follows. If $D$, $D' \in \scd$,  $z \in D$, and $f \in \mathcal{T}(D,D')$, then 
$f \circ \mu_D(z) = \mu_{D'}(f(z))$. This definition is independent of the choice of $f\in \mathcal{T}(D,D')$; indeed, if $f_1$, $f_2 \in \mathcal{T}(D,D')$ with $f_1(z)=f_2(z)=z'$, then $f_1 \circ \mu_D(z) = \mu_{D'}(f_1(z)) = \mu_{D'}(z') =  \mu_{D'}(f_2(z)) = f_2 \circ \mu_D(z)$. The first part of the next proposition follows from a quick change-of-variables, while the second follows immediately from the first as a result of~(\ref{conf-inv-PK}).  See also~\cite[Proposition~5.5]{LawlerSLE}.

\begin{proposition}
\label{nov10prop1}
Suppose that $D$, $D' \in \scd$, and $z \in D$, $w\in \bd D$ with $\bd D$ locally analytic at $w$. If $f \in \mathcal{T}(D,D')$, and $\bd D'$ is locally analytic at $f(w)$, then 
$f \circ \mu_D(z,w) = |f'(w)| \,\mu_{D'}(f(z),f(w))$ and 
$f \circ \mu_D^{\#}(z,w) = \mu_{D'}^{\#}(f(z),f(w))$.
\end{proposition}

Using the interior-to-boundary excursion measure, we now define boundary-to-boundary excursion measure in $\D$, and show that it exists by an explicit calculation. It is then a simple matter to define excursion measure for other simply connected $D$, and to  derive the important conformal invariance formula.

\begin{definition}
If $x$, $y \in \bd \D$, $x \neq y$, then \defines{normalized excursion measure on excursions from $x$ to $y$ in $\D$} is the measure on $\curves$, concentrated on $\curves_x^y(\D)$, defined by
\begin{equation}\label{mar9.eq1}
\lim_{\e \to 0+} \mu_{\D}^{\#}((1-\e)x,y) =:\overline{\mu}_{\bd \D}(x,y),
\end{equation}
where $\mu_{\D}^{\#}(z,y)$ for $z \in \D$, $y \in \bd \D$ is as in~(\ref{normalized-int-bnd-exc-meas}).
\end{definition}

\begin{lemma}
The limit in~(\ref{mar9.eq1}) exists.
\end{lemma}

\begin{proof}
\newcommand{\fa}{f_{\alpha}}
Let $\gamma \in \curves(\D)$ with $\gamma(0)=0$, $\gamma(t_{\gamma})=1$. Let $f_{\alpha}(z)=\frac{z-\alpha}{1-\alpha z}$ for $\alpha \in (-1,1)$ so that $\fa \in \mathcal{T}(\D,\D)$, $\fa(0)=-\alpha$, and both $1$ and $-1$ are fixed points of $\fa$. Using the exact form of the M\"{o}bius transformation $\fa$, a straightforward computation shows that $\lim_{\alpha \to 1} \fa \circ \gamma$ exists in the metric space $(\curves,\metric)$ where $\fa \circ \gamma$ is defined as on page~\pageref{imagepage}. In particular, this shows $\lim_{\e \to 0+} \mu_{\D}^{\#}(-(1-\e),1) =:\overline{\mu}_{\bd \D}(-1,1)$ exists. For other $x$ and $y$, simply use a composition of M\"{o}bius transformations.
\end{proof}

\begin{definition}
We define \defines{excursion measure on excursions from $x$ to $y$ in $\D$} to be the measure on $\curves$, concentrated on $\curves_x^y(\D)$, defined by
$\mu_{\bd \D}(x,y) =  H_{\bd \D}(x,y) \cdot \overline{\mu}_{\bd \D}(x,y)$
where $H_{\bd \D}$ denotes the excursion Poisson kernel.
\end{definition}

Observe that by definition, the mass of excursion measure $\mu_{\bd}(x,y)$
is \emph{defined} to be
$|\mu_{\bd\D}(x,y)| = H_{\bd \D}(x,y)$. Hence, $\mu_{\bd \D}^{\#}(x,y)= \overline{\mu}_{\bd\D}(x,y)$.

\begin{definition}\label{genexcmeas}
Suppose that $D \in \scd$, and  $z$, $w\in \bd D$ with $\bd D$ locally analytic at both $z$ and $w$. 
Let $h \in \mathcal{T}(\D,D)$. \defines{Excursion measure from $z$ to $w$ in $D$} is defined by
\begin{equation}\label{nov11def1}
\mu_{\bd D}(z,w) := \frac{1}{|h'(h^{-1}(z))|\,|h'(h^{-1}(w))|} \,h \circ \mu_{\bd\D}(h^{-1}(z), h^{-1}(w)).
\end{equation}
\end{definition}

A straightforward exercise in the chain rule shows that the definition of $\mu_{\bd D}(z,w)$ given by~(\ref{nov11def1}) does not depend on the choice of $h \in \mathcal{T}(\D,D)$.

\begin{proposition}\label{confinvexcmeas}
Let $D$, $D' \in \scd$, and  let $z$, $w\in \bd D$ with $\bd D$ locally analytic at both $z$, $w$. If $f \in \mathcal{T}(D,D')$, and $\bd D'$ is locally analytic at both $f(z)$, $f(w)$, then
\begin{equation}\label{confinvexcmeas.eq}
f \circ \mu_{\bd D}(z,w) = |f'(z)| \, |f'(w)| \, \mu_{\bd D'}(f(z),f(w))
\end{equation}
and
\begin{equation}\label{confinvexcmeas.eq2}
f \circ \mu_{\bd D}^{\#}(z,w) =  \mu_{\bd D}^{\#}(f(z),f(w)).
\end{equation}
If $f_1$, $f_2 \in \mathcal{T}(D,D')$, then 
$f_1 \circ \mu_{\bd D}(z,w) = f_2 \circ \mu_{\bd D}(z,w)$ so~(\ref{confinvexcmeas.eq}) and~(\ref{confinvexcmeas.eq2}) are independent of the choice of map. In particular, 
\begin{equation*}
\mu_{\bd D}(z,w) = \lim_{\e \to 0+} \frac{1}{\e}\mu_{D}(z+\e\n_z,w).
\end{equation*}
\end{proposition}

\begin{definition}\label{excmeasD}
Suppose that $D \in \scd$. \defines{Excursion measure in $D$} is defined by
$$\mu_{\bd D} := \int_{\bd D}\int_{\bd D} \mu_{\bd D}(z,w) \; |\d w| \; |\d z|.$$
\end{definition}

The conformal invariance of excursion measure is immediate; see~\cite[Proposition~5.8]{LawlerSLE}.

\begin{proposition}[Conformal Invariance]\label{confinvexcmeasD.nowprop}
If $D$, $D' \in \scd$ and $f \in \mathcal{T}(D,D')$, then
$f \circ \mu_{\bd D} = \mu_{\bd D'}$.
\end{proposition}

In fact, it should be noted that we can define excursion measure $\mu_{\bd D}$ for \emph{any} simply connected subset of $\C$ by conformal invariance. The only reason to restrict to $D \in \scd^*$ in the next definition is so that excursions $\gamma \in \curves(D)$ will necessarily have $t_{\gamma} < \infty$. (We will not be concerned with excursions with  $t_{\gamma} = \infty$ in this paper.)

\begin{definition}\label{excmeasD*}
Suppose that $D \in \scd^*$ and $f \in \mathcal{T}(\D, D)$. \defines{Excursion measure in $D$} is defined by $\mu_{\bd D} := f \circ \mu_{\bd \D}$. Furthermore, if $\Gamma$, $\Upsilon \subset \bd D$ with $\overline{\Gamma} \cap \overline{\Upsilon} \neq \emptyset$, define $\mu_{\bd D} (\Gamma, \Upsilon)$ to be the measure $\mu_{\bd D}$ restricted to those excursions $\gamma \in \curves_{\Gamma}^{\Upsilon}(D)$, and define the  \defines{excursion Poisson kernel} $H_{\bd D} (\Gamma, \Upsilon)$ to be its mass; that is, $H_{\bd D} (\Gamma, \Upsilon) := |\mu_{\bd D} (\Gamma, \Upsilon)|$.
\end{definition}

An immediate consequence of these definitions is the following.

\begin{proposition}[Conformal Invariance]\label{excmeasGUprop}
If $D$, $D' \in \scd^*$; $f \in \mathcal{T}(D,D')$; $\Gamma$, $\Upsilon \subset \bd D$ with $\overline{\Gamma} \cap \overline{\Upsilon} \neq \emptyset$; and $\Gamma'$, $\Upsilon'$ are the images under $f$ of $\Gamma$, $\Upsilon$, respectively, then
$f \circ \mu_{\bd D} (\Gamma, \Upsilon) =  \mu_{\bd D'} (\Gamma', \Upsilon')$
and
$H_{\bd D} (\Gamma, \Upsilon) =  H_{\bd D'} (\Gamma', \Upsilon')$. 
\end{proposition}

\subsection{Discrete excursions and discrete excursion measure}

Throughout this section, suppose that $A \in \A$; $w$, $z \in A$; $x$, $y \in \bd A$; and $\G$, $\U \subset \bd A$ with $\overline{\G} \cap \overline{\U} = \emptyset$. Our goal is to define a discrete excursion and formulate the discrete analogues of the previous sections.
If $S_j$ is a simple random walk with $S_0=w$, denote the one-step transition probability $p(w,z) :=  \P^w\{S_1 = z\}$, 
 and define the \defines{discrete Poisson kernel} to be
$h_A(w,y) := \P^{w}\{S_{\tau_A} = y\}$ where $\tau_A := \min \{j>0: S_j \not\in A \}$. As in~\cite[\S 3.1]{Dur1}, 
\begin{equation}\label{discrete_h1}
q(w,z;y) := \P^w\{S_1=z|S_{\tau_A} =y\} = p(w,z)\frac{h_A(z,y)}{h_A(w,y)}
\end{equation}
defines the one-step transition probabilities of simple random walk conditioned to exit $A$ at $y$.
Note that $h_A$ is discrete harmonic, and that~(\ref{discrete_h1}) an example of a discrete $h$-transform.

\begin{definition}
A \defines{discrete excursion in $A$} is a path $\omega := [\omega_0, \omega_1, \ldots, \omega_k]$ where $\omega_0 \in \bd A$, $\omega_k \in \bd A$, $|\omega_i-\omega_{i-1}|=1$ for $i=1, \ldots, k$, and $\omega_i \in A$ for $i=1, \ldots, k-1$, where $2 \le k <\infty$. If $\omega = [\omega_0, \omega_1, \ldots, \omega_k]$, define the \defines{length} of $\omega$, written $|\omega|$, to be $k$.  If $\omega$ is a discrete excursion in $A$ with $\omega_0 \in \Gamma$ and $\omega_k \in \Upsilon$, then $w$ is called  a \defines{$(\Gamma, \Upsilon)$-discrete excursion in $A$}. In particular, if $\omega_0=x$ and $\omega_k=y$, then $\omega$ is called a \defines{discrete excursion from $x$ to $y$ in $A$}.
\end{definition}

Discrete excursions can be generated by starting a simple random walk at $x\in \bd A$, conditioning it to take its first step into $A$, and stopping it at $\tau_A$. Let  the \defines{discrete excursion Poisson kernel} $h_{\bd A}(x,y)$ be given by $h_{\bd A}(x,y) := \P^{x}\{S_{\tau_A} = y , S_1 \in A \}$, and define \defines{discrete excursion measure} to be the measure that assigns weight $4^{-|\omega|}$ to each discrete excursion $\omega$. Denote this measure by $\mu^{\RW}_{\bd A}(\cdot)$ so that $\mu^{\RW}_{\bd A}(\omega) := {4}^{-|\omega|}$. Write $\mu^{\RW}_{\bd A}(x,y)$ to denote the measure on discrete excursions from $x$ to $y$ in $A$, and  $\mu^{\RW}_{\bd A}(\Gamma,\Upsilon) := \sum_{x \in \G} \sum_{y\in \U} \mu^{\RW}_{\bd A}(x,y)$ to denote the measure on $(\Gamma, \Upsilon)$-discrete excursions in $A$. In~\cite{LawW1}, Lawler and Werner defined $\mu^{\RW}_{\bd A}(\omega) := (2\pi\,4^{|\omega|})^{-1}$; this difference only affects things up to a constant.

We want both discrete excursion measure and Brownian excursion measure to be measures on the metric space $(\curves,\metric)$. Consequently, we need to associate to each discrete excursion $\omega$ a curve $\tilde{\omega} \in \curves$.  Suppose that $\omega$ is a discrete excursion in $A$, and let $\cl(A) :=A \cup \bd A$ with associated domain $\widetilde{\cl(A)} \subset \C$. We associate to $\omega$ a curve $\tilde{\omega} \in \curves\left(\widetilde{\cl(A)}\right)$ by setting $t_{\tilde{\omega}} := 2|\omega|$, and 
\newcommand{\intpt}[1]{\lfloor #1 \rfloor}
\begin{equation}\label{discretetoBMexc}
\tilde{\omega}(t) := \omega_{\intpt{t/2}} + \frac{1}{2}(t-\intpt{t})(\omega_{\intpt{t/2}+1}-\omega_{\intpt{t/2}}), \quad 0 \le t \le t_{\tilde{\omega}}.
\end{equation}
In other words, we join the lattice points in order with line segments parallel to the coordinate axes in $\Z^2$, with each segment taking time 2 to traverse. Note that $\tilde{\omega}(0)=\omega_0$ and $\tilde{\omega} (t_{\tilde{\omega}})=\omega_{|\omega|}$.  Using this identification,
if $\omega$ is an excursion from $x$ to $y$ in $A$, then $\mu^{\RW}_{\bd A}(x,y) \in \finmeas(\curves)$ and  $\mu^{\RW}_{\bd A,x,y}(\tilde{\omega}) = 4^{-t_{\tilde{\omega}}}$.  In order to prove discrete excursion measure converges to Brownian excursion measure in the scaling limit, we will consider scaling excursions as the mesh of the lattice becomes finer. See~(\ref{defn-Phi}) in Section~\ref{scale-sect}.

As a consequence of the so-called KMT approximation (see Section~\ref{KMTsection}), it follows that $|B_t - S_{2t}| = O(\log t)$. Complete details may be found in~\cite{KozL} and~\cite{LawlerJose}. Thus, it is simply a matter of \ae sthetics that a random walk path of $|\omega|$ steps take time $2|\omega|$ to traverse: if $\gamma$ is Brownian curve and $\tilde{\omega}$ is as above with $\gamma(0)=\tilde{\omega}(0)$, then $|\gamma(t) - \tilde\omega(t)| = O(\log t)$.

\begin{definition}
Suppose that $A \in \A$ and $x$, $y \in \bd A$. \defines{Discrete excursion measure} $\mu^{\RW}_{\bd A}(x,y)$ is defined to be the measure on $(\curves, \metric)$, concentrated on $V = V(x,y;A) := \{ \gamma \in \curves : \metric(\gamma, \tilde{\omega})=0$ for some discrete excursion $\omega$ from $x$ to $y$ in $A\}$ given by
$\mu^{\RW}_{\bd A}(x,y)(\gamma) := {4}^{-t_{\gamma}}$
for $\gamma \in V$. Note that $\mu^{\RW}_{\bd A}(x,y)(V) = h_{\bd A}(x,y)$.
\end{definition}

\section{The main convergence arguments}\label{Sect4}

\subsection{Carath\'eodory convergence}\label{cara_sect}

\begin{definition}
Fix $r>0$. Suppose that $D_n$ is a sequence of domains with $D_n \in \Dr{r}$ for each $n$. The \defines{kernel} of $D_n$, written $\ker(\{D_n\})$, is the largest domain $D$ containing the origin and having the property that each compact subset of $D$ lies in all but a finite number of the domains $D_n$. Suppose that $\ker(\{D_n\}) =D$. The sequence $D_n$ \defines{converges in the Carath\'eodory sense} to $D$, written $D_n \cara D$, if every subsequence $D_{n_j}$ of $D_n$ has $\ker(\{D_{n_j}\})=D$.
\end{definition}

Recall that a sequence of functions $f_n$ on a domain $D$ converges to a function $f$ \defines{uniformly on compacta of $D$} if for each compact $K \subset D$, $f_n \to f$ uniformly on $K$. The following theorem, which roughly states that convergence  of domains in the Carath\'eodory sense is equivalent to the uniform convergence on compacta of the appropriate Riemann maps, will be our main tool.  A proof may be found in~\cite[Theorem 3.1]{Duren}.  

\begin{theorem}[Carath\'eodory Convergence]\label{caraconvthm}
Suppose that $D_n$ is a sequence of domains with $D_n \in \scd^*$ for each $n$, and let $f_n \in \mathcal{T}(\D,D_n)$ with $f_n(0)=0$, $f_n'(0)>0$. Suppose further that $D \in \scd$ and $f \in \mathcal{T}(\D,D)$ with $f(0)=0$, $f'(0)>0$. Then $f_n \to f$ uniformly on compacta of $\D$ if and only if $D_n \cara D$. 
\end{theorem}

\begin{lemma}\label{CCTcor.nowlem}
Suppose that $D_n \cara D$ with $D_n$, $D \in \scd^*$.  Suppose further that there exists an $E \in \scd^*$ with $D_n \subset E$ for all $n$, and $D \subseteq E$.  If $F:E \to \D$ is the conformal transformation with $F(0)=0$, $F'(0)>0$, then $F(D_n) \cara F(D)$.
\end{lemma}

\begin{proof}
Let $f_n : \D \to D_n$ and let $f : \D \to D$ be conformal transformations mapping 0 to 0 with positive derivative at the origin.  By Theorem~\ref{caraconvthm}, the convergence of $D_n$ to $D$ is equivalent to the uniform convergence of $f_n$ to $f$ on compacta of $\D$. Set $h_n := F \circ f_n$ and $h := F \circ f$, and let $K$ be a compact subset of $\D$. If $z \in K$, then $|h_n(z) - h(z)| = |F(f_n(z))-F(f(z))| \to 0$ uniformly as $n \to \infty$.
\end{proof}

\begin{lemma}\label{jan31lem1.new.nowprop.nowlem}
Suppose that $F_n$, $F$ are conformal mappings of $\D$. Let $D:=F(\D)$. If $F_n \to F$ uniformly on compacta of $\D$, then $F_n \circ F^{-1} \to I$ uniformly on compacta of $D$, where $I: D \to D$ is the identity map $I(z)=z$.
\end{lemma}

\begin{proof}
Let $K' \subset D$ be compact. Let $\eps >0$ be given. Let $K=F^{-1}(K') \subset \D$ which is clearly compact. By uniform convergence, there exists $N=N(\eps,K)$ such that $|F_n(x)-F(x)| < \eps$ for all $n > N$, $x \in K$. If $y \in K'$, then $y=F(x)$ for some $x \in K$.  Hence, if $n > N$, then $|F_n \circ F^{-1}(y)-I(y)| = |F_n(x)-F(x)| < \eps$, and the proof is complete.
\end{proof}

\begin{lemma}\label{limitation-prop.nowlem}
Suppose that $F_n$, $F$ are conformal mappings of the unit disk $\D$, and that $F_n \to F$ uniformly on compacta of $\D$. If $\gamma \in \curves(\D)$ with $|\gamma(0)|<1$ and $|\gamma(t_{\gamma})|<1$, then $\metric(F_n \circ \gamma, F\circ \gamma) \to 0$ as $n \to \infty$.
\end{lemma}

\begin{proof}
Suppose that $\gamma \in \curves(\D)$ with $|\gamma(0)|<1$ and $|\gamma(t_{\gamma})|<1$. Note that $\gamma$ is \emph{not} an excursion in $\D$. Therefore, there necessarily exists a compact set $K \subset \D$ such that $\gamma \in \curves(K)$. Consider $t_{F_n\circ\gamma} = A^n_{t_{\gamma}} =\int_0^{t_{\gamma}} |F_n'(\gamma(r))|^2 \, \d r$ and $t_{F\circ\gamma} = A_{t_{\gamma}} =\int_0^{t_{\gamma}} |F'(\gamma(r))|^2 \, \d r$. Since 
$F_n \to F$ uniformly on compacta of $\D$, we necessarily have that $F_n \to F$ uniformly on $K$. Hence, it follows that $t_{F_n\circ\gamma} \to t_{F\circ\gamma}$. Furthermore,
\begin{align*}
&\sup_{0 \le s \le 1}|F\circ \gamma(t_{F\circ\gamma} s) -F_n \circ \gamma(t_{F_n\circ\gamma} s) |\notag\\
&\;\;\le \sup_{0 \le s \le 1}|F\circ \gamma(t_{F\circ\gamma} s) -F \circ \gamma(t_{F_n\circ\gamma} s) |+|F\circ \gamma(t_{F_n\circ\gamma} s) -F_n \circ \gamma(t_{F_n\circ\gamma} s) |\to 0.
\end{align*}
Taken together, these imply the result.
\end{proof}

\subsection{Construction of approximate domains $\tilde{D}_N$}\label{def_DN_section}

Suppose that $D \in \scd^*$ with $\inrad(D) =1$, and let 
$$D''_N := \left\{ x \in \frac{1}{N}\Z^2 \cap D : \frac{1}{N} \Square_x \subset D \right\},$$
where $\Square_x := x +\left(\,[-1/2,1/2]\times[-1/2,1/2]\,\right)$ is the unit square about $x$. Let $D'_N$ be the connected component of $D''_N$ containing the origin, and set $D_N := D'_N \setminus \bd_i D'_N$. Take $\tilde{D}_N \subset \C$ to be the interior of the union of the scaled squares centred at those $x \in D_N$. We call $D_N$ the \defines{$1/N$-scale discrete approximation to $D$} (with respect to the origin), and we informally refer to $\tilde D_N$ as the associated ``union of squares'' domain; that is, 
\begin{equation}\label{tildeDNdefn}
\tilde{D}_N = \interior \left( \bigcup_{x \in D_N}\frac{1}{N} \, \Square_x \right)\; \text{ and } \; \cl(\tilde{D}_N) := \tilde{D}_N \cup \bd \tilde{D}_N = \bigcup_{x \in {D_N}} \frac{1}{N} \, \Square_x.
\end{equation}
Let $f \in \mathcal{T}(\D,D)$ with $f(0)=0$, $f'(0)>0$. Let $\GD$, $\UD \subset \bd \D$ be (open) boundary arcs with $\overline{\GD} \cap \overline{\UD} = \emptyset$; that is, $\GD := \{e^{i\theta} : \theta_1 < \theta < \theta_2 \}$ and $\UD := \{e^{i\theta} : \theta_3 < \theta < \theta_4 \}$, for some $0 \le \theta_1 < \theta_2 < \theta_3  < \theta_4 < \theta_1+2\pi$. Define $\G \subset \bd D$ to be the image of $\GD$ under $f$, and similarly, let $\U \subset \bd D$ be the image of $\UD$ under $f$. Let $s:=\sep(\G,\U)$ as in Definition~\ref{defnsep.new}, and let $N$ be sufficiently large so that $s \ge \eps_n := n^{-1/48}\log^{2/3} n$ if $n \ge N$. If $f_N \in \mathcal{T}(\D,\tilde{D}_N)$ with $f_N(0)=0$, $f'_N(0)>0$, then define $\tGN$ to be the image of $\GD$ under $f_N$, with $\tUN$ defined similarly. In Theorem~\ref{domain_conv_thm}, we prove $f_N \to f$ uniformly on compacta of $\D$ showing $\tilde{D}_N \cara D$.

 We now define our approximating discrete boundary arcs. If $\tGN \subset \bd \tilde{D}_N$, then associate to $\tGN$ the set $\G_N \subset \bd D_N$ as follows. Let $\G'_N := \{x \in \bd_i D_N : \frac{1}{N}\,\Square_x \cap \tGN \neq \emptyset\}$, and then take
$$\G_N := \left\{ y \in \bd D_N : (x,y) \in \bd_e D_N \text{ with } x\in \G'_N \text{ and } \frac{1}{N}\,\ell_{x,y} \subset \tGN \right\}.$$
Similarly, let $\U_N$ be the discrete boundary arc associated to $\tUN$. Our notation is summarized in the following table.

\begin{center}
\renewcommand{\arraystretch}{1.4}
\begin{tabular}{|c|c|c|c|}\hline
$\D \subset \C$ &$D \subset \C$, $D\in \scd^*$ &$\tilde{D}_N \subset \C$, $\tilde{D}_N \in \scd$  &$D_N \subset \frac{1}{N}\,\Z^2$, $2N D_N \in \A^N$\\\hline
$\Gamma_{\D}, \Upsilon_{\D} \subset \bd \D$   &$\Gamma, \Upsilon \subset \bd D$   &$\tGN, \tUN \subset \bd \tilde{D}_N$ &$\Gamma_N, \Upsilon_N  \subset \bd D_N$    \\\hline
\end{tabular}
\end{center}

Note that by conformal invariance, it is equivalent to specify either $\G$, $\U \subset \bd D$, or $\GD$, $\UD \subset \bd \D$. We have (arbitrarily) chosen the latter.

\subsection{Convergence of domains $\tilde{D}_N$ to $D$}\label{cara-conv-section}

Suppose that $D \in \scd^*$ with $\inrad(D) =1$ and let $D_N$ be the $1/N$-scale discrete approximation to $D$ with associated ``union of squares'' domain $\tilde{D}_N$ as in Section~\ref{def_DN_section}. The following lemmas are an immediate consequence of those definitions.

\begin{lemma}\label{domain_cov_cor1a}
For each $N$, $\tilde{D}_N \in \scd$ with $\cl(\tilde{D}_N) \subset D$. That is, $\tilde{D}_N$ is a simply connected proper subset of $D$ with piecewise analytic boundary. Furthermore, the lattice $\cl(D_N) := D_N \cup \bd D_N \subset D$.
\end{lemma}

\begin{lemma}\label{domain_cov_cor2a}
Suppose that $x \in \bd_i D_N$, $y \in \bd D_N$, and $z \in \bd \tilde{D}_N$. Then $\dist(x,\bd D) \le c_1 \, N^{-1}$, $\dist(y,\bd D) \le c_2 \, N^{-1}$, and $\dist(z,\bd D) \le c_3 \, N^{-1}$ where $c_1 = 2\sqrt{2} + 1/\sqrt{2}$, $c_2 = \sqrt{2} + 1/\sqrt{2}$, and $c_3 = 2\sqrt{2}$.
\end{lemma}

The next proposition follows from the Beurling estimate see~\cite[Proposition~3.79]{LawlerSLE}\label{SLE5}.

\begin{proposition}\label{bdcor2}
If $x \in \bd_i D_N$ and $f \in \mathcal{T}(D,\D)$ with $f(0)=0$, $f'(0)>0$, then
there exists a constant $C$ such that 
$\dist(f(x), \bd \D) \le CN^{-1/2}$  and $f(\tilde{D}_N) \supseteq \{ |z| \le 1 - CN^{-1/2}\}$.
\end{proposition}

We will now establish Theorem~\ref{cara-thm-intro} with the proof of the following result.

\begin{theorem}\label{domain_conv_thm}
The sets $\tilde{D}_N$ as defined by~(\ref{tildeDNdefn}) converge to $D$ in the Carath\'{e}odory sense.
\end{theorem}

The proof of this theorem requires two lemmas.  The first is a simple power series estimate, while the second gives good bounds on the difference of the image of a point under two different maps: the identity map from $\D$ to $\D$, and a map which is ``almost the identity.''

\begin{lemma}\label{log_est}
If $0 \le |z| \le 1/2$, then $|\log(1+z)-z| \le |z|/2$.
\end{lemma}

\begin{proof}
Since 
$$\log(1+z) = \sum_{n=1}^{\infty} \,(-1)^{n-1} \, \frac{1}{n} \, z^n,$$
we have
$$|\log(1+z) -z| \le \sum_{n=2}^{\infty}\, \frac{1}{n} \, |z|^n \le \frac{1}{2}\, |z|\, \sum_{n=1}^{\infty}\, |z|^n \le \frac{1}{2} \, |z|$$
provided that $0 \le |z| \le 1/2$.
\end{proof}

\begin{lemma}\label{lemma_bound}
For $N> 4C^2$, where $C$ is the constant in Lemma~\ref{bdcor2}, suppose that $E_N$ is a domain with $\{|z|\le 1 - CN^{-1/2} \} \subseteq E_N \subseteq \{|z|\le 1 + CN^{-1/2} \}$.  Let $h_N : \D \to E_N$ be the conformal transformation with $h_N(0)=0$ and $h_N'(0) > 0$. Then, there exists a constant $C'$ such that $|h_N(z) - z| \le C' N^{-1/2} \log N$ for $|z|\le 1 - CN^{-1/2}$.
\end{lemma}

\begin{proof}
Without loss of generality, assume that $h_N$ may be extended to an analytic function in a neighbourhood of $\overline{\D}$.  For if this is not the case, we may approximate $h_N$ by $h_{N,r}(z) := r^{-1}h_N(rz)$ and take the limit as $r \to 1-$. From the Schwarz lemma~\cite[page~135]{Ahl1}, we can immediately see that $1-CN^{-1/2} \le h_N'(0) \le 1+CN^{-1/2}$.  Let $\kappa_N(z) := \log[h_N(z)/z]$ so that $\kappa_N = u_N + iv_N$ is analytic on $\D$ with $|u_N(z)| \le (3/2)CN^{-1/2}$ for $|z|=1$ using the estimate from Lemma~\ref{log_est}. Thus, the maximum principle for harmonic functions tells us that $|u_N(z)| \le (3/2)CN^{-1/2}$ for all $|z|\le 1$. We therefore conclude that the partial derivatives of $u_N$ at $z$ are bounded by an absolute constant times $N^{-1/2}\dist(z,\bd \D)^{-1}$; whence $|\kappa_N'(z)| \le C_1 N^{-1/2}(1-|z|)^{-1}$. Writing 
\begin{align*}
\left|\log\left[1+ \frac{h_N(z)-z}{z}\right]\right| = |\kappa_N(z)|
= \left|\kappa_N(0) + \int_0^z \kappa'_N(w) \, \d w\right| 
&\le \frac{C_2}{\sqrt{ N}} \left[1 + \log \frac{1}{1-|z|}\right]
\end{align*}
with $C_2=\max\{C,C_1\}$, we see that if $\epsilon >0$ is such that 
\begin{equation}\label{eqdom1}
\left|\frac{h_N(z)-z}{z}\right| \le \frac{1}{2} \;\; \text{ for } \;\; |z| \le \epsilon,
\end{equation}
then
\begin{equation}\label{eqdom2}
\left|\frac{h_N(z)-z}{z}\right| \le 2 \left|\log\left[1+ \frac{h_N(z)-z}{z}\right] \right| \le 2C_2N^{-1/2}\left[1 + \log \frac{1}{1-|z|}\right].
\end{equation}
Since (\ref{eqdom1}) holds for some $\epsilon>0$, we can iterate (\ref{eqdom2}) to see that (\ref{eqdom2}) must hold for all $|z|$ such that the right side of (\ref{eqdom2}) is less than 1/2.  For $N$ sufficiently large, this includes all $|z| \le 1 - CN^{-1/2}$.
\end{proof}

\begin{proof}[Proof of Theorem~\ref{domain_conv_thm}]
Suppose that $f : D \to \D$ is the conformal transformation with $f(0)=0$, $f'(0)>0$, and let
$\tilde{f}_N : f(\tilde{D}_N) \to \D$ be the conformal transformation with $\tilde{f}_N(0) = 0$, $\tilde{f}'_N(0) > 0$. Let  $F_N : \D \to \tilde{D}_N$ and $F : \D \to D$ be the conformal transformations with $F_N(0)=0$, $F_N'(0)>0$, and $F(0)=0$, $F'(0)>0$, respectively, which are defined by setting $F_N := (\tilde{f}_N \circ f)^{-1}$ and $F :=f^{-1} = (I \circ f)^{-1}$ where $I(z)=z$ is the identity map from $\D$ to $\D$. Finally, let $z \in \D$, and let $w:=\tilde{f}_N^{-1}(z)$ so that $F_N(z)=F(w)$. 

We prove that $\tilde{D}_N \cara D$ by applying Theorem~\ref{caraconvthm} which states that it is sufficient to show $F_N \to F$ uniformly on each compact subset of $\D$. Equivalently, we will show that for each $\delta >0$ sufficiently small, $F_N \to F$ uniformly for all $|z| \le 1- \delta$.
Fix $0 < \delta <1/2$ and choose $M$ so that $M > (3C'\delta^{-1})^{3}$ where $C'$ is the constant in Lemma~\ref{lemma_bound}.  Let $N > M$.  Then by Lemma~\ref{lemma_bound}, we have that for $|z| \le 1 - \delta$, 
$$|w-z| \le \frac{C'\log N}{\sqrt{N}}\, |z| \le \left(\frac{C'\log N}{\sqrt{N}}\cdot \frac{1-\delta}{\delta}\right) \, \delta.$$
Our choice of $M$ guarantees that $C'\delta^{-1}(1-\delta) N^{-1/2}\log N  < 1$ for $N>M$. By~\cite[Corollary~3.25]{LawlerSLE}\label{SLE6}, if for some $0<r<1$, $|w-z| \le r \dist(z,\bd \D)$, then  $$|F(w) - F(z)| \le \frac{4\dist(F(z), \bd D)}{1-r^2} \, |w-z|.$$
Hence, we conclude
$$|F_N(z) - F(z)| = |F(w) - F(z)| \le \left( \frac{4RC'(1-\delta)}{1-\left( \frac{C'\log N}{\sqrt{N}}\cdot\frac{1-\delta}{\delta} \right)^2} \right) \cdot \frac{\log N}{\sqrt{N}}$$
where $R:=\rad(D)$ so that $F_N \to F$ uniformly; whence $\tilde{D}_N \cara D$.
\end{proof}

\begin{corollary}
If $F \in \mathcal{T}(D, \D)$ with $F(0)=0$, $F'(0)>0$, then $F(\tilde{D}_N) \cara \D$. 
\end{corollary}

\begin{proof}
By Lemma~\ref{domain_cov_cor1a}, $\tilde{D}_N \subset D$, so Lemma~\ref{CCTcor.nowlem} yields the result.
\end{proof}

\subsection{Applying results for $A \in \A^N$ to $D_N$}\label{scale-sect}

Suppose that $D \in \scd^*$ with $\inrad(D)=1$. In this section, we combine our construction of $D_N$ with Theorem~\ref{greentheoremB} and~\cite[Theorem~1.1]{KozL} to restate those results for random walk on $D_N$. The most difficult part of this section is keeping track of the notation. 

We begin by mentioning several scaling relationships that will be needed throughout. If $S_n$ is a random walk on $\Z^2$, then for any $r>0$ there is an associated random walk (which we will also denote by $S_n$) on the lattice $r\Z^2$. In other words, there is a one-to-one correspondence between paths from $x$ to $y$ in $A$ on $\Z^2$, and paths from $rx$ to $ry$ in $rA$ on $r\Z^2$. Hence if $A \subset \Z^2$ and $r>0$, then $G_{rA}(rx,ry) = G_A(x,y)$, where the Green's function on the left side is for random walk on the lattice $r\Z^2$, and the Green's function on the right side is for random walk on $\Z^2$. Similarly, we have $h_{rA}(rx,ry) = h_A(x,y)$ for the discrete Poisson kernel, and $h_{\bd rA}(rx,ry) = h_{\bd A}(x,y)$ for the discrete excursion Poisson kernel. 

The conformal invariance of the Green's function for Brownian motion implies that if $D \in \scd^*$ and $r>0$, then $g_{rD}(rx,ry) = g_D(x,y)$. However, from the conformal covariance of the Poisson kernel (Proposition~\ref{PKconf}) and the excursion Poisson kernel (Proposition~\ref{EPKconf}), it follows that $r H_{rD}(rx,ry) = H_D(x,y)$ and $r^2 H_{\bd rD}(rx,ry) = H_{\bd D}(x,y)$. 

Note that a random walk on $D_N$ is taking steps of size $1/N$. Therefore, let $A_N := 2ND_N$ so that $A_N \in \A^N$, and $\tilde{A}_N := \widetilde{(2ND_N)} = 2N\tilde{D}_N \in \scd$. Hence, $z' \in A_N$ if and only if $z := z'/2N \in D_N$. Suppose $x':=2Nx \in A^N$ with $x \in D_N$ and $y':=2Ny \in A^N$ with $y \in D_N$.
Thus, when the above scaling is applied to $\tilde{A}_N$, we conclude that  
\begin{equation}\label{apr12.eqa}
g_{A_N}(x',y') = g_{2ND_N}(2Nx,2Ny) = g_{D_N}(x,y),
\end{equation}
where $g_{D_N}$ denotes the Green's function for Brownian motion in $\tilde{D}_N$. In particular, if $f_{D_N} \in \mathcal{T}(\tilde{D}_N,\D)$ with $f_{D_N}(0)=0$, $f_{D_N}'(0)>0$ and   $f_{A_N} \in \mathcal{T}(\tilde{A}_N,\D)$ with $f_{A_N}(0)=0$, $f_{A_N}'(0)>0$, then 
since $f_{A_N}(x') = f_{D_N}(x)$ and $g_{A_N}(x') = g_{D_N}(x)$, and since we can write
$f_{A_N}(x') = \exp\{-g_{A_N}(x') + i \theta_{A_N}(x')\}$ and $f_{D_N}(x) =  \exp\{-g_{D_N}(x) + i \theta_{D_N}(x)\}$, 
it follows that 
$\theta_{A_N}(x') = \theta_{2N D_N}(2Nx) = \theta_{D_N}(x)$.
Further, in the random walk case,
$G_{A_N}(x',y') = G_{2ND_N}(2Nx,2Ny) = G_{D_N}(x,y)$,
and for $x' := 2Nx \in \bd A_N$ with $x \in \bd D_N$, we have
$h_{\bd A_N}(x',y') = h_{\bd 2N D_N}(2Nx,2Ny) = h_{\bd D_N}(x,y)$;
similarly, $h_{A_N}(0,x') = h_{D_N}(0,x)$, and $h_{A_N}(0,y') = h_{D_N}(0,y)$. 

For $A_N \in \A^N$, let $A^{*}_N := \{ x' \in A_N : g_{A_N}(x') \ge N^{-1/16} \}$ which is consistent with the usage in~\cite{KozL}. If $x' \in A^*_N$, $y' \in A_N$, then Theorem~\ref{greentheoremB} (in particular, its corollary~\cite[Corollary~3.5]{KozL})
 implies that
\begin{equation}\label{jun11.cor.imply}
G_{A_N}(x',y') = \frac{2}{\pi} \,g_{A_N}(x',y') + k_{y'-x'} + O(N^{-7/24}\log N).
\end{equation}
With the above notation in hand, we are finally able to state the following corollary to~(\ref{jun11.cor.imply}).

\begin{corollary}
Let $x \in D_N$ be such that $x' := 2Nx \in (2ND_N)^{*} = A^*_N$, and let $y \in D_N$ with $y' := 2Ny \in A_N$. Then,
$$G_{D_N}(x,y) = \frac{2}{\pi}\, g_{D_N}(x,y) + k_{y'-x'} + O(N^{-7/24}\log N)$$
where $k_z$ is as in~(\ref{defn-kx}).
\end{corollary}

Note that $k_{y'-x'} \le c N^{-3/2}|x-y|^{-3/2}$. Thus, if $|x-y|\ge N^{-29/36}$, then  $k_{y'-x'} = O(N^{-7/24})$, and we have a refined version of the previous corollary.

\begin{corollary}
If $x \in D_N$ with $x' := 2Nx \in A^*_N$,  $y \in D_N$ with $y' := 2Ny \in A_N$, and $|x-y|\ge N^{-29/36}$, then 
\begin{equation*}
G_{D_N}(x,y) = \frac{2}{\pi} \, g_{D_N}(x,y) + O(N^{-7/24}\log N).
\end{equation*}
\end{corollary}

We also have the following corollary to~\cite[Theorem~1.1]{KozL}.

\begin{corollary}\label{mar1.cor1}
If $D \in \scd^*$ with $\inrad(D)=1$, $D_N$ is the $1/N$-scale discrete approximation to $D$, and 
$x$, $y \in \bd D_N$ with $|\theta_{D_N}(x) - \theta_{D_N}(y)| \geq N^{-1/16}\log^2 N$, then 
$$h_{\bd D_N}(x,y) = \frac{(\pi/2)\; h_{D_N}(0,x)\; h_{D_N}(0,y)}{1-\cos(\theta_{D_N}(x)-\theta_{D_N}(y))}
\; \left[1 + O\left(\frac{\log N}{N^{1/16}|\theta_{D_N}(x)-\theta_{D_N}(y)|} \right) \right].$$
\end{corollary}

We now make several observations regarding excursion measure. Suppose $x$, $y \in \bd \tilde{D}_N$ so that $x' := 2Nx$, $y' := 2Ny \in \bd \tilde{A}_N$ as above. If $f(z) = 2Nz$, then $f \in \mathcal{T}(\tilde{D}_N,\tilde{A}_N)$ with $f(0)=0$ and $f'(z)=2N$ for all $z$. Since excursion measure is conformally covariant/invariant, we are able to conclude that
$\mu_{\bd \tilde{D}_N}(x,y) = 4N^2 \mu_{\bd \tilde{A}_N}(x',y')$ and  $\mu_{\bd \tilde{D}_N} = \mu_{\bd \tilde{A}_N}$.
We know from Donsker's theorem that simple random walk converges in the scaling limit to Brownian motion provided that space and time are scaled appropriately.  In order to prove Theorem~\ref{final-thm}, we will need to apply a similar scaling.  Recall from~(\ref{discretetoBMexc}) that if 
$\omega$ is a discrete excursion then we can associate to it a curve $\tilde{\omega} \in \curves$, and that the Brownian scaling map $\brscale_a$ was defined in Example~\ref{brscale_example}. For $N \in \N$, write $\Phi_N := \brscale_{1/(2N)}$ so that
\begin{equation}\label{defn-Phi}
\Phi_{N}\tilde{\omega}(t) = \frac{1}{2N} \; \tilde{\omega}(4N^2t)\;\; \text{ for }\;\; 0 \le t \le t_{\Phi_{N}\tilde{\omega}} = \frac{t_{\tilde{\omega}}}{4N^2} = \frac{|\omega|}{2N^2}.
\end{equation}

\begin{lemma}\label{scale-lemma}
If $\gamma$, $\gamma' \in \curves$, then 
$$
\frac{1}{4N^2}\,\metric(\gamma, \gamma')
\le 
\metric(\Phi_N \gamma, \Phi_N \gamma')
\le
\frac{1}{2N}\,\metric(\gamma, \gamma').$$
\end{lemma}

\begin{proof} From the definitions of $\metric$ and $\Phi_N$ we conclude that
\begin{align*}
\metric(\Phi_N \gamma, \Phi_N \gamma') 
&= \sup_{0 \le s \le 1} |\Phi_N \gamma(st_{\Phi_N \gamma}) -  \Phi_N \gamma'(st_{\Phi_N \gamma'})|
+|t_{\Phi_N \gamma} -  t_{\Phi_N \gamma'}|\\
&= \sup_{0 \le s \le 1} \left|\frac{1}{2N}\, \gamma(st_{\gamma}) -  \frac{1}{2N}\, \gamma'(st_{\gamma'})\right|
+\frac{1}{4N^2}\,|t_{\gamma} -  t_{\gamma'}|
\end{align*}
so the result follows.
\end{proof}

\begin{definition}
Suppose that $D \in \scd^*$ with $\inrad(D)=1$, $D_N$ is the $1/N$-scale discrete approximation to $D$, and $x$, $y \in \bd D_N$. The \defines{$1/N$-scale discrete excursion measure} $\mu^{\RW}_{\bd D_N}(x,y)$ is defined to be the measure on $(\curves, \metric)$, concentrated on $V_N = V_N(x,y;D) := \{ \gamma \in \curves : \metric(\gamma, \Phi_{N}\tilde{\omega})=0$ for some discrete excursion $\omega$ from $2Nx$ to $2Ny$ in $2ND_N\}$ given by
$\mu^{\RW}_{\bd D_N}(x,y)(\gamma) := 4^{-4N^2 t_{\gamma}}=4^{-|\omega|}$
for $\gamma \in V_N$. Finally, if $\G_N$, $\U_N \subset \bd D_N$ with $\overline{\G_N} \cap \overline{\U_N} = \emptyset$, then
$$\mu^{\RW}_{\bd D_N}(\G_N,\U_N) := \sum_{x \in \G_N} \sum_{y \in \U_N} \mu^{\RW}_{\bd D_N}(x,y).$$
\end{definition}

\subsection{Proof of Theorem~\ref{final-thm}}\label{mainresultsection}

In the present section, we establish the following theorem which, as noted in the introduction, may be regarded as the precise formulation of Theorem~\ref{final-thm}.

\begin{theorem}\label{d-star-to-d}
Suppose $D \in \scd^*$ with $\inrad(D)=1$,  and let $\Gamma$, $\Upsilon \subset \bd D$ be open boundary arcs with $\overline{\G} \cap \overline{\U} = \emptyset$. For every $\eps >0$, there exists an $N$ such that
\begin{enumerate}
\item[\emph{\textbf{(a)}}] \hspace{1in}$\displaystyle \left|\;h_{\bd D_N}(\G_N, \U_N) - \frac{1}{4}H_{\bd D}(\G,\U)\;\right| \le \eps$, 
\item[\emph{\textbf{(b)}}] \hspace{1in}$\displaystyle \prohorov\left(\; \mu^{\#}_{\bd \tilde{D}_N}(\tGN, \tUN),\; \mu^{\#}_{\bd D}(\Gamma, \Upsilon)\,\right) \le \eps$, and
\item[\emph{\textbf{(c)}}] \hspace{1in}$\displaystyle \prohorov\left(\;\mu^{\RW,\#}_{\bd D_N}(\G_N, \U_N), \;\mu^{\#}_{\bd \tilde{D}_N}(\tGN, \tUN) \,\right) \le \eps$,
\end{enumerate}
where $D_N$ is the $1/N$-scale discrete approximation to $D$, $\tilde{D}_N \in \scd$ is the ``union of squares'' domain associated to $D_N$, and
 $\G_N$, $\U_N \subset D_N$ are the corresponding discrete boundary arcs with associated boundary arcs $\tGN$, $\tUN \subset \bd \tilde{D}_N$, respectively. In particular, \emph{\textbf{(a)}}, \emph{\textbf{(b)}}, and \emph{\textbf{(c)}} imply that
\begin{equation*}
\lim_{N\to \infty} \prohorov(\;4\,\mu^{\RW}_{\bd D_N}(\G_N, \U_N), \; \mu_{\bd D}(\G,\U) \,) = 0.
\end{equation*}
\end{theorem}

Each of the three parts of Theorem~\ref{d-star-to-d} will be proved in a separate section: in Section~\ref{hHsection} we prove Theorem~\ref{theorem-hH} establishing \textbf{(a)}, in Section~\ref{tailsection} we prove Theorem~\ref{thm-c} establishing \textbf{(b)}, and finally  in Section~\ref{musection} we prove Theorem~\ref{thm-b} establishing \textbf{(c)}. 

\subsubsection{Review of strong approximation of Brownian motion and random walk}\label{KMTsection}

In order to establish Theorem~\ref{d-star-to-d}, it will be necessary to use a strong approximation result which follows from the theorem of Koml\'os, Major, and Tusn\'ady~\cite{KomMT1, KomMT2}. Because of its central r\^ole in the proof, we include the statement here for the convenience of the reader. In what follows, $S_t$ is defined for non-integer $t$ by linear interpolation. 

\begin{theorem}[Koml\'os-Major-Tusn\'ady]\label{KMTthm}
There exists     $c <\infty$  and   a probability
 space $(\Omega, \F, \P)$ on which are    defined a two-dimensional Brownian 
motion $B$ and a two-dimensional simple random walk $S$ with $B_0=S_0$, such
 that for all $\lambda >0$ and each $n \in \N$,
\[\P \left\{\max_{0\leq t   \le n} \left|\frac{1}{\sqrt{2}}\,B_t - S_t\right| >
 c \,(\lambda + 1) \, \log n  \right\} < c \, n^{-\lambda}.\]
\end{theorem}

The proofs of the following two results may be found in~\cite[Corollary~3.2]{KozL} and~\cite[Proposition~3.3]{KozL}, respectively.

\begin{corollary}\label{strapprox}
There  exist $C< \infty$
and  a probability space $(\Omega, \F, \P)$ on which are  
defined a two-dimensional Brownian motion $B$ and a
 two-dimensional simple random walk $S$ with $B_0=S_0$  such that
$$\P \left\{\max_{0 \le t \le \sigma_n} 
  \left|\frac 1 {\sqrt 2}B_t - S_t\right| > C \log n \right\} = O(n^{-10}),$$
where $\sigma^{1}_n := \inf\{t : |S_t-S_0|
 \ge n^8\}$, $\sigma^{2}_n := \inf\{t : |B_t-B_0| \ge n^8\}$,
 and $\sigma_n := \sigma^{1}_n \vee \sigma^{2}_n$.
\end{corollary}

\begin{proposition}[Strong Approximation]\label{sept23.thm1.nowprop}
There exists a constant $c$ such that for every $n$, 
a Brownian motion $B$ and a simple random walk $S$ can be defined on
the same probability space so that if $A \in \A^n$, $1 < r \leq n^{20}$, 
and $x \in A$ with $|x| \leq n^3$, then
$\P^x\{|B_{T_A} - S_{\tau_A}| \geq c r \log n \} \leq  c  r^{-1/2}$.
\end{proposition}

By combining the strong approximation with Theorem~\ref{KMTthm}, the following estimate is easily deduced.

\begin{proposition}\label{theconjecture}
There exists a decreasing sequence $\dn \downarrow 0$ such that if $A \in \A^n$ with associated ``union of squares'' domain $\tilde{A} \in \scd$, and $\G$, $\U \subset \bd A$ with $\overline{\G} \cap \overline{\U} = \emptyset$ and associated boundary arcs $\tilde{\G}$, $\tilde{\U} \subset \bd \tilde{A}$, then $h_A(0,\Gamma) = H_{\tilde{A}}(0,\tilde{\Gamma}) + O(\dn)$, and $h_A(0,\Upsilon) = H_{\tilde{A}}(0,\tilde{\Upsilon}) + O(\dn)$. Consequently,
$h_A(0,\Gamma)\,h_A(0,\Upsilon) =H_{\tilde{A}}(0,\tilde{\Gamma})\,H_{\tilde{A}}(0,\tilde{\Upsilon}) + O(\dn)$
where the error term depends on both $\Gamma$, $\Upsilon$.
\end{proposition}

\begin{proof}
If $c$ is the constant in Proposition~\ref{sept23.thm1.nowprop}, and $V$ is the set $V := \{ x \in \bd A : \dist(x,\G) \le cn^{1/8}\log n \}$, then $h_A(0,\G) = H_{\tilde{A}}(0,\tilde{V}) + O(n^{-1/16})$. However, a simple gambler's ruin estimate for Brownian motion shows that $H_{\tilde{A}}(0,\tilde{V}) =  H_{\tilde{A}}(0,\tilde{\G}) + O(n^{-7/8}\log n)$, so  the result follows with $\dn = n^{-7/8}\log n$. 
\end{proof}

\subsubsection{Convergence of $4 h_{\bd D_N}(\G_N,\U_N)$ to $H_{\bd D}(\G,\U)$}\label{hHsection}

The goal of the present section is to prove that if $D \in \scd^*$ with $\inrad(D)=1$, and $\G$, $\U \subset \bd D$ are disjoint open boundary arcs, then $4 h_{\bd D_N}(\G_N,\U_N) \to H_{\bd D}(\G,\U)$ using the notation from Section~\ref{def_DN_section}, therefore establishing Theorem~\ref{d-star-to-d} \textbf{(a)}.  

It follows from the exact form of the excursion Poisson kernel in $\D$~\cite[Example~5.6]{LawlerSLE} that if $D \in \scd$ and $x$, $y \in \bd D$ with $\bd D$ locally analytic at $x$ and $y$, then 
\begin{equation}\label{EPK-PKform-eqn-new}
H_{\bd D}(x,y) =\frac{2\pi \, H_{D}(0,x) \, H_{D}(0,y)}{1- \cos(\theta_D(x)-\theta_D(y))}.
\end{equation}
For further details, see also~\cite[Example~2.14]{KozL}.

\begin{lemma}\label{EPK-PKform-general} 
If $D \in \scd^*$ and $\G$, $\U \subset \bd D$ with $\overline{\G} \cap \overline{\U} \neq \emptyset$, 
then
$$\frac{2\pi H_{D}(0,\G)\,H_{D}(0,\U)}{1-\cos(\spr(\Gamma,\Upsilon))} \le H_{\bd D}(\G,\U)
\le \frac{2\pi H_{D}(0,\G)\,H_{D}(0,\U)}{1-\cos(\sep(\Gamma,\Upsilon))}$$
where $H_{\bd D}(\Gamma, \Upsilon)$ is as in Definition~\ref{excmeasD*}, and $H_D(0,\Gamma)$, $H_D(0,\Upsilon)$ are as in~(\ref{PKjune05}).
\end{lemma}

\begin{proof}
Suppose first that $D \in \scd$, and that $\G$, $\U$ are analytic open boundary arcs. 
Then from~(\ref{EPK-PKform-eqn-new}), we conclude that for all $x \in \G$ and for all $y \in \U$, 
$$\frac{2\pi H_{D}(0,x) H_{D}(0,y)}{1-\cos(\spr(\Gamma,\Upsilon))} 
\le H_{\bd D}(x,y)  
\le \frac{2\pi H_{D}(0,x) H_{D}(0,y)}{1-\cos(\sep(\Gamma,\Upsilon))}.$$
Since $D \in \scd$, Proposition~\ref{EPKconf} implies that 
$$\frac{2\pi H_{D}(0,\G)\,H_{D}(0,\U)} {1-\cos(\spr(\Gamma,\Upsilon))} \le H_{\bd D}(\G,\U)
\le \frac{2\pi H_{D}(0,\G)\,H_{D}(0,\U)}{1-\cos(\sep(\Gamma,\Upsilon))}.$$
Now, suppose $D' \in \scd^*$, and let $f \in \mathcal{T}(D,D')$. Write $\G'$, $\U'$ for the images under $f$ of $\G$, $\U$, respectively. Conformal invariance yields $H_{D}(0,\G) = H_{D'}(0,\G')$ and $H_{D}(0,\U) = H_{D'}(0,\U')$. (Indeed this holds for all domains $D\in \scd^*$ since $\bd D$ is regular.) From Proposition~\ref{EPKconf}, it follows that $H_{\bd D}(\G,\U) = H_{\bd D'}(\G',\U')$; whence the proof is complete.
\end{proof}

Let $f \in \mathcal{T}(\D,D)$ with $f(0)=0$, $f'(0)>0$. Analogous to Section~\ref{def_DN_section}, by rotating\footnote{Both the excursion Poisson kernel for $\D$ and excursion measure in $\D$ are rotationally invariant.} if necessary, it is possible to find $0 \le \theta_1 < \theta_2 < \theta_3 < \theta_4 < 2\pi$ such that $\G$, $\U$, are the images under $f$ of $\GD$, $\UD$, respectively, where $\GD := \{e^{i\theta}:  \theta_1 < \theta < \theta_2\}$ and $\UD := \{e^{i\theta'}:  \theta_3 < \theta' < \theta_4\}$. Define the \defines{length of $\G$}, written $\lenD{\Gamma}$, to be length of $\GD$ so that $\lenD{\Gamma} := \theta_2-\theta_1$. Similarly define $\lenD{\Upsilon} := \theta_4-\theta_3$. Note that our notion of length is simply harmonic measure so that while $\Gamma$ may not even be rectifiable, $\lenD{\Gamma}$ always exists. An easy estimate shows that  if $(\theta_3-\theta_2)$, $(\theta_4-\theta_1)$ are fixed, then
$$\frac{1-\cos(\theta_3-\theta_2)}{1-\cos(\theta_4-\theta_1)} = 1 +O(\theta_4-\theta_3)+ O(\theta_2-\theta_1)$$
as $(\theta_4-\theta_3)\to 0$, $(\theta_2-\theta_1)\to 0$, and hence,
as $\lenD{\Upsilon} \to 0$, $\lenD{\Gamma} \to 0$,
\begin{equation}\label{mar3.eq1}
\frac{1-\cos(\sep(\Gamma,\Upsilon))}{1-\cos(\spr(\Gamma,\Upsilon))} =  1 +O(\lenD{\Upsilon})+ O(\lenD{\Gamma}).
\end{equation}
Thus, we have proved the following lemma. 

\begin{lemma}\label{mar3.lem1}
If $D \in \scd^*$, then for any $\eta >0$ there exist open boundary arcs $\G$, $\U \subset \bd D$ with $\overline{\G} \cap \overline{\U}=\emptyset$ such that 
$$1 \le \frac{1-\cos(\spr(\G,\U))}{1-\cos(\sep(\G,\U))} \le 1+\eta.$$
\end{lemma}

Note that the lower bound  holds automatically by the definitions of separation and spread.

\begin{theorem}\label{theorem-hH}
For every $D \in \scd^*$ with $\inrad(D)=1$, and for every pair of open boundary arcs $\G$, $\U \subset \bd D$ with $\overline{\G} \cap \overline{\U} \neq \emptyset$, 
if $D_N$ is the $1/N$-scale discrete approximation to $D$, and $\G_{N}$, $\U_{N}$ are the discrete approximations to $\G$, $\U$, respectively, as in Section~\ref{def_DN_section}, then $4 h_{\bd D_N}(\G_{N},\U_{N})\to H_{\bd D}(\G,\U)$.
\end{theorem}

\begin{proof}
Consider $D \in \scd^*$, and let $\G$, $\U \subset \bd D$ be (open) boundary arcs with $\overline{\G} \cap \overline{\U} \neq \emptyset$. Find $M$ so that $\sep(\G,\U) \ge \eps_N:= N^{-1/48}\log^{2/3} N$ for $N \ge M$. Throughout this section, let $N \ge M$. Let $D_N$ be the $1/N$-scale discrete approximation to $D$ with associated ``union of squares'' domain $\tilde{D}_N$, and let $\tGN$, $\tUN \subset \bd \tilde{D}_N$ with associated discrete boundary arcs $\G_N$, $\U_N \subset \bd D_N$. From the definitions of separation and spread, and from Corollary~\ref{mar1.cor1},  since $\G$ and $\U$ are fixed so that $\sep(\G,\U) = O(1)$, it follows that
$$
\frac{(\pi/2) h_{D_N}(0,x) h_{D_N}(0,y)}{1-\cos(\spr(\Gamma,\Upsilon))} [1 + O(\e_N^3)]
\le h_{\bd D_N}(x,y) \le 
\frac{(\pi/2) h_{D_N}(0,x) h_{D_N}(0,y)}{1-\cos(\sep(\Gamma,\Upsilon))} [1 + O(\e_N^3)].
$$
Summing over all $x \in \G_N$ and all $y \in \U_N$ yields
$$
\frac{h_{D_N}(0,\G_N) h_{D_N}(0,\U_N)}{1-\cos(\spr(\Gamma,\Upsilon))} [1 + O(\e_N^3)]
 \le \frac{2}{\pi}\,  h_{\bd D_N}(\G_N,\U_N)\le 
\frac{h_{D_N}(0,\G_N) h_{D_N}(0,\U_N)}{1-\cos(\sep(\Gamma,\Upsilon))} [1 + O(\e_N^3)].
$$
where we write 
$h_{\bd D_N}(\G_N,\U_N) :=\sum_{x \in \G_N} \sum_{y \in \U_N} h_{\bd D_N}(x,y)$
and similarly for $h_{D_N}(0,\G_N)$ and $h_{D_N}(0,\U_N)$. However, from Proposition~\ref{theconjecture},
$$h_{D_N}(0,\G_N)\, h_{D_N}(0,\U_N) = H_{\tilde{D}_N}(0,\tGN)\,H_{\tilde{D}_N}(0,\tUN) + O(\dN),$$
where $\dN := N^{-7/8}\log N$, so that we conclude
\begin{align*}
&\left[\frac{H_{\tilde{D}_N}(0,\tGN)\,H_{\tilde{D}_N}(0,\tUN)}{1-\cos(\spr(\Gamma,\Upsilon))}  +O(\dN)\right] \cdot [1 + O(\e_N^3)]\\
&\qquad \le \frac{2}{\pi} \,h_{\bd D_N}(\G_N,\U_N)\le 
\left[\frac{H_{\tilde{D}_N}(0,\tGN)\,H_{\tilde{D}_N}(0,\tUN)}{1-\cos(\sep(\Gamma,\Upsilon))} +O(\dN) \right]\cdot [1 + O(\e_N^3)].
\end{align*}
Now, as we let $N \to \infty$, it follows that
\begin{align*}
&\frac{H_{D}(0,\G)\,H_{D}(0,\U)}{1-\cos(\spr(\Gamma,\Upsilon))} \le \frac{2}{\pi} \, \liminf_{N\to \infty} h_{\bd D_N}(\G_N,\U_N)\\
&\qquad\le \frac{2}{\pi} \, \limsup_{N \to \infty} h_{\bd D_N}(\G_N,\U_N) \le \frac{H_{D}(0,\G)\,H_{D}(0,\U)}{1-\cos(\sep(\Gamma,\Upsilon))}.
\end{align*}
However, Lemma~\ref{EPK-PKform-general} implies that
\begin{align*}
&\frac{1-\cos(\sep(\Gamma,\Upsilon))}{1-\cos(\spr(\Gamma,\Upsilon))} H_{\bd D}(\G,\U) \le 4 \; \liminf_{N\to \infty} h_{\bd D_N}(\G_N,\U_N)\\
&\qquad \le 4 \; \limsup_{N \to \infty} h_{\bd D_N}(\G_N,\U_N) \le \frac{1-\cos(\spr(\Gamma,\Upsilon))}{1-\cos(\sep(\Gamma,\Upsilon))} H_{\bd D}(\G,\U).
\end{align*}
For any $\eta >0$, let $\{\G_i\}$, $\{\U_j\}$ be finite partitions of $\G$, $\U$, respectively, with
$$1 \le \frac{1-\cos(\spr(\G_i,\U_j))}{1-\cos(\sep(\G_i,\U_j))} \le 1+\eta.$$
Note that such a partitioning is possible by Lemma~\ref{mar3.lem1}. Hence, the equation above becomes
\begin{align*}
&\frac{1-\cos(\sep(\G_i,\U_j))}{1-\cos(\spr(\G_i,\U_j))} H_{\bd D}(\G_i,\U_j)\le 4 \; \liminf_{N\to \infty} h_{\bd D_N}(\G_{N,i},\U_{N,j})\\
&\qquad \le 4 \; \limsup_{N \to \infty} h_{\bd D_N}(\G_{N,i},\U_{N,j}) \le \frac{1-\cos(\spr(\G_i,\U_j))}{1-\cos(\sep(\G_i,\U_j))} H_{\bd D}(\G_i,\U_j).
\end{align*}
Summing over $i$ and $j$ and noting that
$$\sum_{i} \sum_{j} H_{\bd D}(\G_i,\U_j) = H_{\bd D}(\G,\U)\; \text{ and }\; \sum_{i} \sum_{j} h_{\bd D_N}(\G_{N,i},\U_{N,j}) = h_{\bd D_N}(\G_{N},\U_{N})$$
since 
$\{\G_{N,i}\}$, $\{\U_{N,j}\}$ partition $\{\G_N\}$, $\{\U_N\}$, respectively, gives
\begin{align*}
&(1+\eta)^{-1} H_{\bd D}(\G,\U) \le 4 \; \liminf_{N\to \infty} h_{\bd D_N}(\G_{N},\U_{N})\\
&\qquad \le 4 \; \limsup_{N \to \infty} h_{\bd D_N}(\G_{N},\U_{N}) \le (1+\eta) H_{\bd D}(\G,\U).
\end{align*}
Since $\eta>0$ was arbitrary, we conclude $4 h_{\bd D_N}(\G_{N},\U_{N}) \to  H_{\bd D}(\G,\U)$ as $N \to \infty$.
\end{proof}

\subsubsection{Convergence of $\mu^{\#}_{\bd \tilde{D}_N}(\tGN, \tUN)$ to $\mu^{\#}_{\bd D}(\G, \U)$}\label{tailsection}

We now prove Theorem~\ref{d-star-to-d} \textbf{(b)} via a result which basically says that an excursion in $D$ can be thought of as an excursion in $\tilde{D}_N$ with Brownian tails.

\begin{theorem}\label{thm-c}
For every $D \in \scd^*$ with $\inrad(D)=1$, and for every pair of open boundary arcs $\G$, $\U \subset \bd D$ with $\overline{\G} \cap \overline{\U} \neq \emptyset$, 
\begin{equation}\label{thm-c-eq}
\lim_{N \to \infty}\prohorov\left(\,\mu^{\#}_{\bd \tilde{D}_N}(\tGN, \tUN), \,\mu^{\#}_{\bd D}(\G, \U) \,\right) = 0
\end{equation}
where $D_N$ is the $1/N$-scale discrete approximation to $D$ with associated domain $\tilde{D}_N \in \scd$, and corresponding boundary arcs $\tGN$, $\tUN \subset \bd \tilde{D}_N$ as in~Section~\ref{def_DN_section}. 
\end{theorem}

By conformal invariance, we can define excursion measure $\mu^{\#}_{\bd D}(\Gamma, \Upsilon)$ to be either the measure
$f \circ \mu^{\#}_{\bd \D}(\GD, \UD)$ for $f \in \mathcal{T}(\D, D)$, or $\mu_{\bd D}$ restricted to those excursions $\gamma \in \curves_{\Gamma}^{\Upsilon}(D)$ (and normalized by $H_{\bd D}(\Gamma, \Upsilon)$). Also using conformal invariance, we have $\mu^{\#}_{\bd \tilde{D}_N}(\tGN, \tUN) = f_N \circ \mu^{\#}_{\bd \D}(\GD, \UD)$ for $f_N \in \mathcal{T}(\D, \tilde{D}_N)$,
so that we conclude
\begin{equation}\label{jan31eq1}
\mu^{\#}_{\bd \tilde{D}_N}(\tGN, \tUN) =  (f_N \circ f^{-1}) \circ  \mu^{\#}_{\bd D}(\Gamma, \Upsilon).
\end{equation}
In order to show the convergence of the masses $h_{\bd D_N}(\GN,\UN)$ to $H_{\bd D}(\G,\U)$, the intermediate step of showing  
\begin{equation}\label{apr7.eq1}
\lim_{N \to \infty}H_{\bd \tilde{D}_N}(\tGN, \tUN) = H_{\bd D}(\Gamma, \Upsilon)
\end{equation}
is unnecessary as a consequence of the conformal invariance of the excursion Poisson kernel: $H_{\bd \tilde{D}_N}(\tGN, \tUN) = H_{\bd D}(\Gamma, \Upsilon)$. However, in contrast to the excursion Poisson kernel, it is not simply a matter of applying the conformal invariance of excursion measure to conclude that (cf.~Lemma~\ref{limitation-prop.nowlem})
\begin{equation}\label{apr12.eq1}
\prohorov\left(\,(f_N \circ f^{-1}) \circ  \mu^{\#}_{\bd D}(\Gamma, \Upsilon)\, , \, \mu^{\#}_{\bd D}(\G, \U)\, \right) \to 0.
\end{equation}

Suppose that $D \in \scd^*$ with $\inrad(D)=1$, and associated ``union of squares'' domain $\tilde{D}_N$.
As mentioned in Lemma~\ref{domain_cov_cor2a}, if $z \in \bd \tilde{D}_N$, then $\dist(z,\bd D) \le 2\sqrt{2} \, N^{-1}$. It follows from the Beurling estimates (see~\cite[Proposition 3.79]{LawlerSLE} and~\cite[Lemma~5.3]{LawlerJose}) that Brownian motion started at $z$ is likely to exit $D$ quickly and nearby; that is, 
\begin{equation}\label{tail-est}
\P^{z} \{ \diam B[0,T_D] \ge N^{-1/2} \} \le CN^{-1/4} \;\; \text{ and }\;\; \P^{z} \{ T_D \ge N^{-1/2} \} \le CN^{-3/8}.
\end{equation}
Unfortunately, if $z \in \tGN$, it may be extremely unlikely that $\{ B_{T_D} \in \G \}$. This will be the case, for example, if $z$ and $\G$ are on opposite sides of a ``channel'' (or ``fjord''). However, since $\tilde{D}_N \cara D$ by Theorem~\ref{domain_conv_thm}, for \emph{fixed} $D \in \scd^*$, fixed disjoint open boundary arcs $\G$, $\U$, and for every $\eps > 0$, there exists an $N$ such that  
$\max \{\dist(\tGN,\G), \dist(\tUN,\U)\} < \eps$. The following is then a consequence of~(\ref{tail-est}) and easy bounds on the Poisson kernel.

\begin{lemma}\label{apr12.lema}
For every $\eps > 0$, there exists an $N$ such that for all $z \in \tGN$,
\begin{equation}\label{tail-est-2}
\P^{z} \{ T_D \ge \eps \;\text{ or }\; \diam B[0,T_D] \ge \eps \;\text{ or }\; B_{T_D} \not\in \G_{\eps} \} \le \eps
\end{equation}
where $\G_{\eps} := \{z \in \bd D: \dist(z,\G) \le \eps \}$.
\end{lemma}

\begin{proof}[Proof of Theorem~\ref{thm-c}]
Suppose that $\gamma : [0,t_{\gamma}] \to \C$ is a $(\tGN,\tUN)$-excursion in $\tilde{D}_N$. Let $b_2 : [0,t_{b_2}] \to \C$ be a Brownian motion started at $\gamma(t_{\gamma})$ and stopped at $t_{b_2} := \inf\{t : b_2(t) \in D\}$, its hitting time of $\bd D$. Let $b' : [0,t_{b'}] \to \C$ be an independent Brownian motion started at $\gamma(0)$, stopped at $t_{b'} := \inf\{t : b'(t) \in D\}$, and set $b_1(t) :=  b'(t_{b'}-t)$. If $\zeta := b_1 \oplus \gamma \oplus b_2$, then by construction $\zeta : [0,t_{\zeta}] \to \C$ has $\zeta(0)\in \bd D$, $\zeta(t_{\zeta})\in \bd D$, $0<t_{\zeta}<\infty$, and $\zeta(0,t_{\zeta}) \subset D$.  In other words, $\zeta$ is an excursion in $D$. Unfortunately, $\zeta$ is not necessarily a $(\Gamma,\Upsilon)$-excursion in $D$, but with high probability is very close to one. Indeed, if we denote by $\nu_{\bd \tilde{D}_N}(\Gamma,\Upsilon)$ the probability measure on paths obtained by this \emph{$(\tGN, \tUN)$-excursion in $\tilde{D}_N$ plus Brownian tails} procedure, then it follows from~(\ref{tail-est-2}) and Proposition~\ref{d-p-prop} that for every $\eps > 0$ there exists an $N$ such that   
$$\P\{\metric(\zeta,\gamma) \ge \eps\} \le \eps \;\;\text{ and therefore }\;\;\prohorov\left(\nu_{\bd \tilde{D}_N}(\Gamma,\Upsilon), \mu^{\#}_{\bd \tilde{D}_N}(\tGN,\tUN)\right) \le \eps.$$
The proof is completed by noting that 
$\prohorov(\nu_{\bd \tilde{D}_N}(\Gamma,\Upsilon), \mu^{\#}_{\bd D}(\G,\U) ) \to 0\,$
as a consequence of Proposition~\ref{confinvexcmeas}: $(\G,\U)$-Brownian excursions in $D$ are generated by starting $\eps$ from $\G$ inside $D$ and conditioning the Brownian motion to exit $D$ at $\U$. 
\end{proof}

As in the discussion preceding Theorem~\ref{final-thm}, we can use~(\ref{thm-c-eq}) and~(\ref{apr7.eq1}) to define the convergence of the infinite measures $\mu_{\bd \tilde{D}_N}$ to $\mu_{\bd D}$.

\begin{theorem}\label{final-thm2}
If $D \in \scd^*$ with $\inrad(D)=1$, then $\prohorov(\,\mu_{\bd \tilde{D}_N}, \mu_{\bd D} ) \to 0$ where $D_N$ is the $1/N$-scale discrete approximation to $D$ with associated  domain $\tilde{D}_N$.
\end{theorem}

It must be noted, however, that by Proposition~\ref{confinvexcmeasD.nowprop} and Definition~\ref{excmeasD*}, we \emph{define} excursion measure $\mu_{\bd D}$ for $D \in D^*$ by conformal invariance. Let $f_N \in \mathcal{T}(\D,\tilde{D}_N)$ as above, and also suppose that $f \in \mathcal{T}(\D,D)$. Hence, 
$\mu_{\bd \tilde{D}_N} := f_N \circ \mu_{\bd \D}$ and $\mu_{\bd D} := f \circ \mu_{\bd \D}$
so that $\mu_{\bd \tilde{D}_N} = (f_N \circ f^{-1}) \circ \mu_{\bd D}$
as in~(\ref{jan31eq1}). Thus, we can rephrase the conclusion of Theorem~\ref{final-thm2} as 
$\prohorov(\, (f_N \circ f^{-1}) \circ \mu_{\bd D}, \; \mu_{\bd D} \, ) \to 0$; compare this with~(\ref{apr12.eq1}).

\subsubsection{Estimating $\prohorov\left(\mu^{\RW,\#}_{\bd D_N}(\G_N, \U_N), \;\mu^{\#}_{\bd \tilde{D}_N}(\tGN, \tUN) \right)$}\label{musection}

In this section we establish Theorem~\ref{d-star-to-d} \textbf{(c)} by proving the following result.

\begin{theorem}\label{thm-b}
For every $D \in \scd^*$ with $\inrad(D)=1$, for every pair of open boundary arcs $\G$, $\U \subset \bd D$ with $\overline{\G} \cap \overline{\U} \neq \emptyset$, and for every $\eps >0$, there exists an $N$ such that
\begin{equation}\label{apr6.eq1}
\prohorov\left(\mu^{\RW,\#}_{\bd D_N}(\G_N, \U_N), \mu^{\#}_{\bd \tilde{D}_N}(\tGN, \tUN) \right) \le \eps
\end{equation}
where $D_N$ is the $1/N$-scale discrete approximation to $D$ with associated domain $\tilde{D}_N \in \scd$ and corresponding boundary arcs $\GN$, $\UN \subset \bd D_N$; $\tGN$, $\tUN \subset \bd \tilde{D}_N$ as in~Section~\ref{def_DN_section}. 
\end{theorem}

In order to prove~(\ref{apr6.eq1}), it will be necessary to use the strong approximation of Proposition~\ref{sept23.thm1.nowprop}. Hence, let $A_N := 2N D_N$ so that $A_N \in \A^N$, and write $\G_{N,A} :=2N\GN$, $\U_{N,A} :=2N\UN \subset \bd A_N$ for the corresponding boundary arcs. Suppose further that $N$ is chosen large enough so that $\dist(\G_{N,A},\U_{N,A}) \ge N^{15/16}$. Since $D \in \scd^*$, it follows that $A_N$ is necessarily bounded so that $\rad(A_N) \asymp \inrad(A_N) \asymp N$, and furthermore, $|\U_{N,A}| \asymp |\G_{N,A}| \asymp N$ where all of the constants may depend on $D$. 

\begin{proof}[Proof of Theorem~\ref{thm-b}]
Suppose that $x \in A_N^* :=\{ x \in A_N : g_{A_N}(x) \ge N^{-1/16}\}$, and let $S$ be a simple random walk with $S_0=x$. 
As in the proof of~\cite[Corollary~3.5]{KozL}, it follows from the Beurling estimate that $\dist(x,\bd A) \ge CN^{7/8}$.
Hence, a straightforward gambler's ruin estimate shows that $\P^x \{ S_{\tau} \in \U_{N,A} \} \asymp N^{-1/16}$ where $\tau = \tau_{A_N} := \min \{ j : S_j \in \bd A\}$. The coupling of Brownian motion and random walk provided by Corollary~\ref{strapprox} is so strong that even conditioning on the rare event $\{S_{\tau} \in \U_{N,A} \}$ does not uncouple the processes. Hence, there exists a Brownian motion $B$, a simple random walk $S$ with $B_0=S_0=x$,  and a constant $C$ such that
\begin{equation}\label{apr6.eq2}
\P^{x} \left\{ \sup_{0 \le t \le \tau} \left|\frac{1}{\sqrt{2}}\,B_t - S_{t}\right| \ge C \log N  \,\Big|\, S_{\tau} \in \U_{N,A} \right\} \le C N^{-8}.
\end{equation}
The strong approximation (Proposition~\ref{sept23.thm1.nowprop}) allows us to conclude that conditioned on the event $\{S_{\tau} \in \U_{N,A} \}$, Brownian motion and simple random walk starting $N^{7/8}$ away from the boundary still exit near each other; that is,
\begin{equation}\label{apr6.eq3}
\P^{x} \left\{  |B_{T} - S_{\tau}| \ge C N^{1/4}\log N \,|\, S_{\tau} \in \U_{N,A} \right\} \le C N^{-1/16}
\end{equation}
where $T = T_{A_N}:= \inf \{t : B_t \in \bd \tilde{A}_N \}$. The time version of the Beurling estimate~\cite[Lemma~5.3]{LawlerJose} says that 
$\P^{x} \{  |T - \tau| \ge  r^{2}\dist(x,\bd \tilde{A})^{2} \} \le C r^{-1/2}$. Hence,
\begin{equation}\label{apr6.eq4}
\P^{x} \left\{  |T - \tau| \ge C N^{1/2}\log^2 N \,|\, S_{\tau} \in \U_{N,A} \right\} \le C N^{-1/16}.
\end{equation}
We can now use Proposition~\ref{d-p-prop} to deduce statements about convergence in $\prohorov$ from statements about convergence in $\metric$. In particular, let $\gamma : [0,t_{\gamma}] \to \C$ be given by $t_{\gamma} := T$, $\gamma(t) := B_t$, $0 \le t \le t_{\gamma}$, and associate to the random walk $S$ the curve $\tilde{\w}:[0,t_{\tilde{\w}}] \to \C$ as in~(\ref{discretetoBMexc}), so that from~(\ref{apr6.eq2}),~(\ref{apr6.eq3}), and~(\ref{apr6.eq4}), we conclude that
$\P \{ \metric(\gamma,\tilde{\w}) \ge C N^{1/2}\log^2 N \} \le CN^{-1/16}$,
and using Lemma~\ref{scale-lemma}, we can scale our results to $D_N$:
\begin{equation}\label{apr6.eqa}
\P \left\{ \metric(\Phi_N\gamma,\Phi_N\tilde{\w}) \ge C N^{-1/2}\log^2 N \right\} \le  \P \left\{ \metric(\gamma,\tilde{\w}) \ge C N^{1/2}\log^2 N \right\} \le N^{-1/16}
\end{equation}
where $\Phi_N := \brscale_{1/(2N)}$ is the Brownian scaling map as in~(\ref{defn-Phi}). Let $V_{N,A}$ be the set $V_{N,A} := \{x \in \bd A_N: \dist(x,\U_{N,A}) \le  CN^{1/4}\log N \}$, let $\tilde{V}_{N,A}$ be the associated subset of $\bd \tilde{A}_N$, and let $2N\tilde{V}_N = \tilde{V}_{N,A}$. It then follows that $\mathcal{L}(\Phi_N\tilde{\w}) = \mu_{D_N}^{\RW,\#}(x,\U_{N})$ and $\mathcal{L}(\Phi_N\gamma) = \mu_{\tilde{D}_N}^{\#}(x,\tilde{V}_{N})$. Since $N^{-1/2}\log N \ll N^{-1/16}$, Proposition~\ref{d-p-prop} and~(\ref{apr6.eqa}) yield
\begin{equation}\label{apr6.eq5}
\prohorov\left(\mu_{D_N}^{\RW,\#}(x,\U_{N}), \mu_{\tilde{D}_N}^{\#}(x,\tilde{V}_{N}) \right) \le CN^{-1/16}.
\end{equation}
As in the proof of Proposition~\ref{theconjecture},  $H_{\tilde{D}_N}(x,\tilde{V}_{N}) = H_{\tilde{D}_N}(x,\tilde{\U}_{N}) + O(N^{-3/4}\log N)$, so it follows that
\begin{equation}\label{apr6.eq6}
\prohorov\left(\mu_{\tilde{D}_N}^{\#}(x,\tilde{\U}_{N}), \mu_{\tilde{D}_N}^{\#}(x,\tilde{V}_{N}) \right) \le CN^{-3/4}\log N.
\end{equation}
Combining~(\ref{apr6.eq5}) and~(\ref{apr6.eq6}) then yields
$\prohorov(\mu_{D_N}^{\RW,\#}(x,\U_{N}),\mu_{\tilde{D}_N}^{\#}(x,\tilde{\U}_{N})) \le  CN^{-1/16}$,
and, in particular, if $y \in \A_N^*$ with $|x-y| \le C\log N$, then 
\begin{equation}\label{apr9.kmteq1}
\prohorov \left(\mu_{D_N}^{\RW,\#}(x,\U_{N}),\mu_{\tilde{D}_N}^{\#}(y,\tilde{\U}_{N}) \right) \le  CN^{-1/16}.
\end{equation}
To complete the proof, suppose that $S'$ is a simple random walk on the scaled lattice $\frac{1}{2N}\,\Z^2$, and 
let $D_N^{*} := \frac{1}{2N}\,A_N^*$ so that $D_N^* = \{ z \in D_N : g_{D_N}(z) \ge N^{-1/16}\}$ by~(\ref{apr12.eqa}) where $g_{D_N}$ is the Green's function for Brownian motion on $\tilde{D}_N$. Also recall from Theorem~\ref{domain_conv_thm} that $\tilde{D}_N \cara D$. Hence, if $\eta_N = \eta(D,N) := \min \{j \ge 0: S'_j \in D_N^* \cup  D_N^c \}$  and $x \in D_N\setminus D_N^*$, then 
it follows from~\cite[Lemma~3.11]{KozL} that for every $\eps > 0$, there exists an $N$ such that 
\begin{equation}\label{apr12.eq-1}
\P^x \left\{ \eta_N \ge \eps \,\big|\, S'_{\eta_N} \in D_N^* \right\} \le \eps.
\end{equation}
Furthermore, using~\cite[Lemma~3.11]{KozL} again, we can find constants $C$, $\alpha$ such that
\begin{equation}\label{apr12.eq-2}
\P^x\left\{\max_{0 \leq j \leq \eta-1}|f_{D_N}(S'_j)-f_{D_N}(x)|\geq  N^{-1/16} \log N\right\} \leq C \:N^{-\alpha},
\end{equation}
and
\begin{equation}\label{apr12.eq-3}
\P^x \left\{ \,|f_{D_N}(S'_\eta)-f_{D_N}(x)|\geq N^{-1/16}\log N \,\big|\, S'_\eta \in D_N^{*} \right\} \leq C \: N^{-\alpha}.
\end{equation}
Suppose further that $\tilde{B}$ is a Brownian motion started at $x \in D_N \setminus D_N^*$. As in Lemma~\ref{apr12.lema},  
if $\tilde{\eta}_N = \tilde{\eta}(D,N) := \inf \{t \ge 0: \tilde{B}_t \in \tilde{D}_N^* \cup  \tilde{D}_N^c \}$,
then for every $\eps >0$, there exists an $N$ such that
\begin{equation}\label{apr12.eq-4}
\P^{x} \left\{ \tilde{\eta}_N \ge \eps \;\text{ or }\; \diam B[0,\tilde{\eta_N}] \ge \eps \,\big|\, B_{\tilde{\eta}_N} \in \tilde{D}_N^* \right\} \le \eps.
\end{equation}
If we let $\tilde{\gamma} : [0, t_{\tilde{\gamma}}] \to \C$ be given by $t_{\tilde{\gamma}} := \tilde{\eta}_N$, $\tilde{\gamma}(t) := \tilde{B}_t$, $0 \le t \le t_{\tilde{\gamma}}$, and associate to the (scaled) random walk $S'$ the (scaled) curve $\tilde{\w}' : [0, t_{\tilde{\w}'}] \to \C$ as in~(\ref{defn-Phi}) (i.e., Brownian scaled in both time and space), then letting $\underset{\widetilde{}}{\gamma} := \tilde{\gamma} \oplus \Phi_N\gamma$ and $\underset{\widetilde{}}{\w} := \tilde{\w}' \oplus \Phi_N\tilde{\w}$ we see that $\mathcal{L}(\underset{\widetilde{}}{\gamma}) = \mu_{\bd \tilde{D}_N}^{\#}(\tGN,\tUN)$
and $\mathcal{L}(\underset{\widetilde{}}{\w}) = \mu_{\bd D_N}^{\RW, \#}(\G_N,\U_N)$. Hence,
by combining~(\ref{apr12.eq-1}),~(\ref{apr12.eq-2}),~(\ref{apr12.eq-3}), and~(\ref{apr12.eq-4}) with~(\ref{apr9.kmteq1}),
we conclude that for every $\eps>0$, there exists an $N$ with
\begin{equation*}
\prohorov\left(\mu_{\bd D_N}^{\RW,\#}(\G_{N},\U_{N}),\mu_{\bd \tilde{D}_N}^{\#}(\tilde{\G}_N,\tilde{\U}_{N})\right) \le  \eps.\qedhere
\end{equation*}
\end{proof} 

\section*{Acknowledgements}

Much of this research was done by the author in his Ph.D.~dissertation~\cite{KozThesis} under the supervision of Greg Lawler.  The author wishes to express his gratitude to Prof.~Lawler for his continued guidance and support. Thanks are also due to Christophe Garban for valuable comments, to the Banff International Research Station for Mathematical Innovation and Discovery where the final writing of this paper was done, and to Christian G.~Bene\v{s} and Jos\'{e} A.~Trujillo Ferreras for many fruitful discussions.


\end{document}